\documentclass[11pt]{article}

\usepackage{fancyhdr}
\usepackage{amsfonts}
\usepackage{amsmath}
\usepackage{amssymb}
\usepackage{amsthm}
\usepackage{tikz}

\usepackage{amscd}
\usepackage{amstext}

\usepackage{oldgerm} 

\usepackage{fancyhdr}
\usepackage{amsfonts}
\usepackage{amsmath}
\usepackage{amssymb}
\usepackage{amsthm}

\usepackage{amscd}
\usepackage{amstext}
\usepackage{color}

\usepackage{graphics}
\usepackage{graphicx}

\newlength{\fixboxwidth}
\newcommand{\N}{\mathbb{N}}

\newcommand{\com}{\mathbb{C}}

\newcommand{\cs}{{\mathcal S}}

\newcommand{\cf}{{\mathcal F}}
\newcommand{\cfi}{{\cf}^{-1}}

\newcommand{\supp}{{\rm supp \, }}
\newcommand{\sgn}{{\rm sgn \, }}

\newcommand{\be}{\begin{equation}}
\newcommand{\ee}{\end{equation}}
\newcommand{\beq}{\begin{eqnarray}}
\newcommand{\beqq}{\begin{eqnarray*}}
\newcommand{\eeq}{\end{eqnarray}}
\newcommand{\eeqq}{\end{eqnarray*}}

\newtheorem{satz}{Theorem}

\newtheorem{rem}{Remark}
\newtheorem{defi}{Definition}
\newtheorem{lem}{Lemma}
 
\newtheorem{prop}{Proposition}

\begin{document}


\title{Quarklet Characterizations for Triebel-Lizorkin Spaces}
\author{Marc Hovemann$^{1,}$\thanks{The work of this author has been supported by Deutsche Forschungsgemeinschaft (DFG), grant $DA 360/24-1$.} $ \ $ and Stephan Dahlke$^{2}$ }
\date{\today}

\maketitle

\begin{center}
{\scriptsize  Institute of Mathematics, Philipps-University Marburg, Hans-Meerwein-Stra{\ss}e 6, 35043 Marburg, Germany \\
Email: 1: hovemann@mathematik.uni-marburg.de   $  \qquad $     2:   dahlke@mathematik.uni-marburg.de  }
\end{center}

\textbf{Abstract.} In this paper we prove that under some conditions on the parameters the one-dimensional Triebel-Lizorkin spaces $ F^{s}_{r,q}(\mathbb{R})   $ can be described in terms of quarklets. So for functions from Triebel-Lizorkin spaces we obtain a quarkonial decomposition as well as a new equivalent quasi-norm. For that purpose we use quarklets that are constructed out of biorthogonal compactly supported Cohen-Daubechies-Feauveau spline wavelets, where the primal generator is a cardinal B-spline. Moreover we introduce some sequence spaces apposite to our quarklet system and study their properties. Finally we also obtain a quarklet characterization for the Triebel-Lizorkin-Morrey spaces $ \mathcal{E}^{s}_{u,r,q}(\mathbb{R})   $.

\vspace{0.2 cm }


\vspace{0.2 cm }


\section{Introduction and main results}

Nowadays the Triebel-Lizorkin spaces $F^s_{r,q} (\mathbb{R})$ are a well-established tool to describe the regularity of functions and distributions. They have been introduced around 1970 by Lizorkin (see \cite{Liz1} and \cite{Liz2}) and Triebel (see \cite{Tr73}). Later these function spaces have been investigated in detail in the books of Triebel, see for example \cite{Tr83}, \cite{Tr92}, \cite{Tr06} and \cite{Tr20}. Usually the Triebel-Lizorkin spaces are defined by using Fourier analytical tools. Aside from this very important are their decompositions in terms of atoms and smooth wavelets. For that we refer to the Chapters 1.5 and 3.1 in \cite{Tr06} as well as to \cite{Tr04a} and Chapter 7 in \cite{FJ2}. In the last years many more decomposition techniques for the Triebel-Lizorkin spaces have been developed. In particular the concept of subatomic and quarkonial decomposition seems to be a quite powerful approach. In Chapter I.2 in \cite{Tr01} these decompositions are based on a partition of unity whose elements are not only translated and dilated but also multiplied with polynomials up to a certain order. In the present paper we work with quarklets that are constructed out of biorthogonal compactly supported Cohen-Daubechies-Feauveau spline wavelets, where the primal generator is a cardinal B-spline. For the theory of such biorthogonal wavelets we refer to Chapter 6.A in \cite{CoDau}. Roughly speaking our quarklets $ \psi_{p}  $ are a linear combination of translated cardinal B-splines that are multiplied with some monomial of degree $ p \in \mathbb{N}_{0}  $. For $ j \in \mathbb{N}  $ and $ k \in \mathbb{Z}   $ by $  \psi_{p,j,k}  $ we denote the dilated and translated versions of our quarklets. Precise definitions can be found in the Definitions \ref{Bquark} and \ref{def_quarklet} below. One important motivation to study quarkonial decompositions are some
very interesting applications in the context of the numerical treatment
of elliptic partial differential equations. Here finite-element-methods (FEM) are well-established tools. The classical $h$-FEM relies on a space refinement, whereas for $p$-methods the polynomial degree of the ansatz functions varies. A combination of both methods is also possible and called $hp$-FEM. For an overview of FEM we refer to \cite{Ci}, \cite{Ha} and \cite{Sch}. When we think on wavelets there exist adaptive wavelet methods that are guaranteed to converge with optimal order, see for example \cite{CaDaDe}. Those strategies can be interpreted as $h$-methods. Therefore the question arises whether
it is possible to design $hp$-versions of adaptive wavelet schemes. This issue directly leads to our quarklets. Some first results concerning adaptive quarklet approximation can be found in \cite{DaFKRaa}, \cite{DaKRaa}, \cite{DaORaa} and \cite{DaRaaS}. In this paper it is our main goal to use quarklets to describe one-dimensional Triebel-Lizorkin spaces $F^{s}_{r,q}(\mathbb{R})  $ via equivalent quasi-norms. In connection with that we also find quarklet representations for functions out of those function spaces. With other words we prove that under some conditions on the parameters the quarklet system is stable in $ F^{s}_{r,q}(\mathbb{R})   $. This is clearly an important advantage in numerical applications,
in particular since usually the regularity of the solution in
Triebel-Lizorkin spaces describes the approximation order that can be
achieved by adaptive numerical schemes. Details concerning this topic can be found in \cite{DaHaSchSi}. With that in mind we are now prepared to formulate our first main result. It tells us that under some conditions on the parameters the one-dimensional Triebel-Lizorkin spaces $ F^{s}_{r,q}(\mathbb{R})   $ can be characterized in terms of quarklets. 
 
\begin{satz}\label{mainresult1}
Let $  s \in \mathbb{R}  $, $  0 < r < \infty $, $ 0 < q < \infty  $ and $ m \in \mathbb{N}  $ with $ m \geq 2  $. Moreover we assume that the parameters fulfill one of the following conditions.
\begin{itemize}
\item[(I)] We have
\begin{align*}
\max \Big(0, \frac 1r -1 , \frac 1q -1  \Big) < s  < m - 1 .
\end{align*}

\item[(II)] We have $ s \geq m - 1  $ with $ r \geq 1   $ and
\begin{align*}
\max \Big(0 , \frac 1q -1  \Big) < s  < \min \Big ( m - 1 + \frac{1}{r} , m - 1 + \frac{1}{q}    \Big ).
\end{align*}
For the fine index $q$ we assume $ \frac{1}{q} < \min \big( m ,  s + 1 - mr(s + 1 - m) \big)  $.

\item[(III)] We have $ s \geq m - 1   $ with $ \frac{1}{m} < r < 1   $ and 
\begin{align*}
\max \Big(  \frac 1r -1  , \frac 1q -1  \Big) < s  < \min \Big ( m - 1 + \frac{r}{q}  , m  \Big ).
\end{align*}
For the fine index we assume $ \frac{1}{q} < \min \big ( m , - sm + m^{2} + \frac{s}{r} - \frac{m}{r} + \frac{1}{r} \big )   $.

\item[(IV)] We have $ 1 < r < \infty   $ and $ 1 < q < \infty    $ as well as $ -m + 1 < s  < 0 $.

\item[(V)] We have $ 1 < r < \infty   $, $ 1 < q < \infty    $ and $ s \leq - m + 1  $ with
\begin{align*}
 \max \Big ( - m + \frac{1}{r} , - m + \frac{1}{q}    \Big ) < s  \leq - m + 1.
\end{align*}
For the fine index $q$ we assume $ \frac{1}{q} > s  + m(1 - \frac{1}{r})(-s + 1 - m)   $.

\end{itemize}
Let $ f \in \mathcal{S}'(\mathbb{R})   $. Then we have   $ f \in F^{s}_{r,q}(\mathbb{R})  $ if and only if $ f $ can be represented as
\begin{equation}\label{rep_main}
f = \sum_{p \geq 0} \sum_{j = -1}^{\infty} \sum_{k \in \mathbb{Z}} c_{p,j,k} \psi_{p,j,k}
\end{equation}
with convergence in $ \mathcal{S}'(\mathbb{R})  $, where we have that
\begin{equation}\label{rep_main_seqnorm}
\Big\| \Big[ \sum_{p = 0}^{\infty}  \sum_{j= -1}^\infty \sum_{k \in \mathbb{Z}}  ( p + 1 )^{\sgn(s) (2m+3)^{2} q}  2^{jsq} 2^{\frac{jq}{2}}    |c_{p,j,k}|^{q}  |\chi_{j,k}(x)|^{q}   \Big]^{\frac{1}{q}}\Big|  L_{r} (\mathbb{R}) \Big\|    
\end{equation}
is finite. Moreover the quasi-norms $ \Vert f \vert F^{s}_{r,q}(\mathbb{R})  \Vert    $ and
\begin{align*}
\inf_{\eqref{rep_main}} \Big\| \Big[ \sum_{p = 0}^{\infty}  \sum_{j= -1}^\infty \sum_{k \in \mathbb{Z}}  ( p + 1 )^{ \sgn(s)  (2m+3)^{2} q}  2^{jsq} 2^{\frac{jq}{2}}    |c_{p,j,k}|^{q}  |\chi_{j,k}(x)|^{q}   \Big]^{\frac{1}{q}}\Big|  L_{r} (\mathbb{R}) \Big\|  
\end{align*}
are equivalent. Here the infimum is taken over all sequences $ \{ c_{p,j,k}  \}_{p \in \mathbb{N}_{0}, j \in \mathbb{N}_{0} \cup \{-1\}, k \in \mathbb{Z} }   $ such that \eqref{rep_main} is fulfilled.
\end{satz}
Here $ \sgn    $ stands for the sign function. Notice that in general the representation in \eqref{rep_main} is not unique. The expression in \eqref{rep_main_seqnorm} is strongly connected with the original sequence spaces associated with the Triebel-Lizorkin spaces. In fact we recover the original sequence spaces when we put $ c_{p,j,k} = 0   $ whenever $ p > 0  $. Let us add some remarks concerning the conditions on the parameters. We observe that in $ (II)$, $(III)$ and $(V)$ some additional conditions concerning the fine index $q$ show up. They appear because of technical reasons, and it might be that they can be weakened by using a refined proof technique. This problem will be studied in a forthcoming paper. More comments referring to the necessity of some conditions on the parameters can be found in Remark \ref{rem_nec1}. In the diagram Figure 1 we illustrate the situation in Theorem \ref{mainresult1} for the special case $r=q$. The numbers $ (I)-(V)  $ refer to those in the theorem. Notice that for $ (II)  $, $(III)$ and $(V)$ not all conditions from the theorem can be found in the diagram.

\begin{tikzpicture}[thick]
\draw[->] (2,0) -- (8,0) ;
\draw[->] (3,-4) -- (3,4) ;
\draw (4,-0.1) -- (4,0.1) ;
\draw (6,-0.1) -- (6,0.1) ;
\draw (2.9,3) -- (3.1,3) ;
\draw (2.9,2) -- (3.1,2) ;
\draw (2.9,-3) -- (3.1,-3) ;
\draw (2.9,-2) -- (3.1,-2) ;

\draw (3,2) -- (4,3) ;
\draw (4,3) -- (6,3) ;
\draw (4,0) -- (6,2) ;
\draw (6,2) -- (6,3) ;

\draw (4,0) -- (4,-2) ;
\draw (4,-2) -- (3,-3) ;




\draw[dotted] (3,2) -- (6,2) ;

\draw[dotted] (4,2) -- (4,3) ;

\draw[dotted] (3,-2) -- (4,-2) ;



\node at (4.2,-0.4) {$1$} ;
\node at (2.3,2) {$m-1$} ;
\node at (2.3,3) {$m$} ;
\node at (2.3,-2) {$-m+1$} ;
\node at (2.3,-3) {$-m$} ;
\node at (6,-0.3) {$m$} ;
\node at (8.3,0) {$\frac{1}{r}$} ;
\node at (3,4.3) {$s$} ;
\node at (5.9,1) {$ s = \frac{1}{r}-1  $} ;
\node at (4,1) {$(I)$} ;
\node at (3.7,2.3) {$(II)$} ;
\node at (4.9,2.3) {$(III)$} ;
\node at (3.5,-1) {$(IV)$} ;
\node at (3.3,-2.2) {$(V)$} ;
\node at (10,- 3.5){Figure 1: The situation in Theorem \ref{mainresult1}.} ;
\node at (0,0){};

\end{tikzpicture}

In recent times an increasing number of authors also deals with the so-called Triebel-Lizorkin-Morrey spaces $ \mathcal{E}^{s}_{u,r,q}(\mathbb{R})    $. They are a generalization of the original Triebel-Lizorkin spaces. So we have $ \mathcal{E}^{s}_{r,r,q}(\mathbb{R}) = F^{s}_{r,q}(\mathbb{R})  $. The Triebel-Lizorkin-Morrey spaces have been introduced by Tang and Xu in 2005, see \cite{TangXu}. For a precise definition and further explanations we refer to Section \ref{sec_mor}. It turns out that also the spaces $ \mathcal{E}^{s}_{u,r,q}(\mathbb{R})  $ can be described in terms of quarklets. So our second main result reads as follows.

\begin{satz}\label{mainresult2}
Let $  s \in \mathbb{R}  $, $  0 < r  \leq u < \infty $, $ 0 < q < \infty  $ and $ m \in \mathbb{N}  $ with $ m \geq 2  $. Moreover we assume
\begin{align*}
\max \Big(0, \frac 1r -1 , \frac 1q -1  \Big) < s  < m - 1 .
\end{align*}
Let $ f \in \mathcal{S}'(\mathbb{R})   $. Then we have   $ f \in \mathcal{E}^{s}_{u,r,q}(\mathbb{R})  $ if and only if $ f $ can be represented as
\begin{equation}\label{rep_main2}
f = \sum_{p \geq 0} \sum_{j = -1}^{\infty} \sum_{k \in \mathbb{Z}} c_{p,j,k} \psi_{p,j,k}
\end{equation}
with convergence in $ \mathcal{S}'(\mathbb{R})  $, where we have
\begin{align*}
\Big\| \Big[ \sum_{p = 0}^{\infty}  \sum_{j= -1}^\infty \sum_{k \in \mathbb{Z}}  ( p + 1 )^{(2m+3)^{2}q}  2^{jsq} 2^{\frac{jq}{2}}    |c_{p,j,k}|^{q}  |\chi_{j,k}(x)|^{q}   \Big]^{\frac{1}{q}}\Big|  \mathcal{M}^{u}_{r}(\mathbb{R}) \Big\|   < \infty .
\end{align*}
Moreover the quasi-norms $ \Vert f \vert \mathcal{E}^{s}_{u,r,q}(\mathbb{R})  \Vert    $ and
\begin{align*}
\inf_{\eqref{rep_main2}} \Big\| \Big[ \sum_{p = 0}^{\infty}  \sum_{j= -1}^\infty \sum_{k \in \mathbb{Z}}  ( p + 1 )^{(2m+3)^{2}q}  2^{jsq} 2^{\frac{jq}{2}}    |c_{p,j,k}|^{q}  |\chi_{j,k}(x)|^{q}   \Big]^{\frac{1}{q}}\Big| \mathcal{M}^{u}_{r}(\mathbb{R})  \Big\|  
\end{align*}
are equivalent. Here the infimum is taken over all sequences $ \{ c_{p,j,k}  \}   $ such that \eqref{rep_main2} is fulfilled.
\end{satz}

When we compare the Theorems \ref{mainresult1} and \ref{mainresult2} we observe that the result for the Triebel-Lizorkin-Morrey spaces only covers a smaller range of parameters. The main reason for this is that some tools we used for the proof of Theorem \ref{mainresult1} seem to be not available for the more difficult Morrey case. 

This paper is organized in the following way. In Section \ref{sec_def} we give some basic definitions. So we define the Triebel-Lizorkin spaces and explain the concept of quarklets. In Section \ref{sec_quas>0} we prove our main result, that is Theorem \ref{mainresult1}, for the case $ s>0  $. For that purpose we work with differences and apply complex interpolation. In Section \ref{sec_quas<0} we prove our main result for $ s<0$. Here we use some duality arguments. In Section \ref{sec_mor} we deal with the Triebel-Lizorkin-Morrey spaces and prove our second main result Theorem \ref{mainresult2}. 

\section*{Notation}

We start by introducing some notation. As usual $\mathbb{N}$ denotes the natural numbers, $\mathbb{N}_0$ the natural numbers including $0$, $\mathbb{Z}$ the integers and $\mathbb{R}$ the real numbers. We put
\[
 B(x,t) := \{y\in \mathbb{R} : \quad |x-y|< t\}\, , \qquad x \in \mathbb{R} \, , \quad t>0.
\]
All functions are assumed to be complex-valued, i.\,e. we consider functions $f:~ \mathbb{R} \to \com$. Let $\mathcal{S}(\mathbb{R})$ be the collection of all Schwartz functions on $\mathbb{R}$ endowed with the usual topology and denote by $\mathcal{S}'(\mathbb{R})$ its topological dual, namely the space of all bounded linear functionals on $\mathcal{S}(\mathbb{R})$ endowed with the weak $\ast$-topology. The symbol $\cf$ refers to  the Fourier transform,
$\cfi$ to its inverse transform, both defined on $\cs'(\mathbb{R})$. Almost all function spaces which we consider in this paper are subspaces of $\cs'(\mathbb{R})$, i.\,e. spaces of equivalence classes with respect to almost everywhere equality. However, if such an equivalence class  contains a continuous representative, then usually we work with this representative and call also the equivalence class a continuous function. 
By $C^\infty_0(\mathbb{R})$ we mean the set of all infinitely often differentiable functions on $\mathbb{R}$ with compact support. For $ 0 < r \leq \infty $ by $ L_{r}(\mathbb{R})  $ we denote the usual Lebesgue spaces. Given a function $ f \in  L_{r}(\mathbb{R}) $ we use the symbol $ \Vert f \vert L_{r}(\mathbb{R}) \Vert   $ for the associated quasi-norm. When we have $ f, g \in L_{2}(\mathbb{R})  $ we use the abbreviation
\begin{align*}
\left\langle f , g     \right\rangle_{L_{2}(\mathbb{R})} = \int_{\mathbb{R}} \overline{f(x)} g(x) dx .
\end{align*}
For $ j, j' \in \mathbb{Z}   $ the symbol  $ \delta_{j , j'}   $ refers to the Kronecker delta. For two quasi-Banach spaces $X$ and $Y$ the operator norm of a linear operator $T:\, X\to Y$ is denoted by $\|T| X \rightarrow  Y \|$. We write $ X \hookrightarrow Y $ if $ X \subset Y $ and the natural embedding of $ X $ into $ Y $ is continuous. For all $r\in(0,\infty )$ and $ q \in (0 , \infty ] $ we write 
\[
\sigma_r:= \,  \max \Big(0, \frac 1r - 1\Big) \qquad \mbox{and}\qquad 
 \sigma_{r,q}:= \,  \max \Big(0, \frac 1r -1 , \frac 1q - 1 \Big) \, .
\]
The symbols  $C, C_1, c, c_{1} \ldots $ denote  positive constants that depend only on the fixed parameters $s,r,q$ and possibly on some additional parameters. Unless otherwise stated their values may vary from line to line. For $ j \in \mathbb{Z}  $ and $ k \in \mathbb{Z}  $ we define the dyadic cubes $ Q_{j,k} $ via $ Q_{j,k} := 2^{-j} ( [0,1) + k   )$. By $ \chi_{j,k}  $ we denote the characteristic function of a cube $ Q_{j,k} $.

\section{Basic definitions: Function spaces and quarklets}\label{sec_def}

\subsection{Besov spaces and Triebel-Lizorkin spaces}

In this section we want to recall the definition of the one-dimensional Triebel-Lizorkin spaces $  F^{s}_{r,q}(\mathbb{R})  $. For that purpose we need a so-called smooth dyadic resolution of the unity. Such a system of functions can be constructed in the following way. Let $\lambda_0 \in C_0^{\infty}(\mathbb{R})$ be a non-negative function such that 
 $\lambda_0(x) = 1$ if $|x|\leq 1$ and $ \lambda_0 (x) =0$ if $|x|\geq \frac{3}{2}$. For $k\in \N$ we define $  \lambda_k(x) := \lambda_0(2^{-k}x)- \lambda_0(2^{-k+1}x) $ for all $ x \in  \mathbb{R} $. Then for every $ x \in  \mathbb{R} $ we observe $ \sum_{k=0}^\infty \lambda_k(x) = 1 $. Moreover for $  k \in \N  $ we find $ \supp \lambda_k \subset \big\{x\in \mathbb{R}  : \: 2^{k-1}\le |x|\le 3 \cdot 2^{k-1}\big\} $. Because of those properties we call the system $(\lambda_k)_{k\in \N_0 }$ a smooth dyadic resolution of the unity on $ \mathbb{R}  $. Clearly by the Paley-Wiener-Schwarz Theorem we find, that $\cfi[\lambda_{k}\, \cf f]$ with $k \in \N_0$ is a smooth function for all $f\in \cs'(\mathbb{R})$. Now we are prepared to define the Triebel-Lizorkin spaces. A definition of the Besov spaces $  B^{s}_{r,q}(\mathbb{R})   $ that are strongly connected with the Triebel-Lizorkin spaces also can be found in what follows.

\begin{defi}\label{def-lt}
Let $(\lambda_k)_{k\in \N_0 }$ be the above system. 
Let $0<  r < \infty$, $ 0 < q \leq \infty$  and $s \in \mathbb{R}$. 
\begin{itemize}
\item[(i)]
Then the Besov space $ B^{s}_{r,q}(\mathbb{R})   $ is the collection of all $  f \in \mathcal{S}'(\mathbb{R})  $ such that
\begin{align*}
\Vert f \vert B^{s}_{r,q}(\mathbb{R})   \Vert :=   \Big ( \sum_{k = 0}^{\infty} 2^{ksq}   \Vert \mathcal{F}^{-1}[\lambda_{k} \mathcal{F}f]  \vert L_{r}(\mathbb{R})   \Vert  ^{q} \Big   )^{\frac{1}{q}}  < \infty .
\end{align*}
In the special case $ q = \infty  $ the sum is replaced by a supremum in the usual way.

\item[(ii)]
The Triebel-Lizorkin space $ F^{s}_{r,q}(\mathbb{R})$ is the collection of all $f \in \mathcal{S}'(\mathbb{R})$ such that
\begin{align*}
 \|\, f \, |F^s_{r,q}(\mathbb{R})\| :=
         \Big\| \Big(\sum\limits_{k=0}^{\infty} 2^{k s q}\, 
         \vert \cfi[\lambda_{k}  \cf f](\, \cdot \, ) \vert^q \Big)^{\frac{1}{q}} \, \Big|L_r(\mathbb{R})\Big\|
          <\infty .
\end{align*}   
In the special case $ q = \infty$ the usual modifications have to be made.
\end{itemize}
\end{defi}
There exists an overwhelming amount of literature concerning the properties of the spaces $ F^{s}_{r,q}(\mathbb{R})   $. Let us refer at least to \cite{Tr83}, \cite{Tr92}, \cite{Tr06} and \cite{Tr20}. For us it is of special interest that the Triebel-Lizorkin spaces also can be described in terms of differences. To give an explanation concerning this topic we have to introduce some additional notation. Let $f:~ \mathbb{R} \to \com$ be a function. Then we define the difference of first order by 
\begin{align*}
\Delta_h^1 f (x) = f(x+h)-f(x)
\end{align*}
for $x,h \in \mathbb{R}$.
Higher order differences are defined via $ \Delta^{N+1}_h f (x) =  \Delta_h^1 (\Delta^{N}_h f) (x) $ with $ x,h \in  \mathbb{R} $ and $ N \in \N\, $. Such differences can be used to formulate equivalent characterizations for the spaces $ F^s_{r,q}(\mathbb{R})  $. Let us recall the following result.

\begin{lem}\label{discrete}
Let $0 < r < \infty$, $0 < q \le \infty$, $ 0 < w \leq \infty   $, $N \in \N$ and
\begin{align*}
\max \Big(0, \frac 1r -1 , \frac 1q - 1 , \frac 1r - \frac{1}{w} , \frac 1q - \frac{1}{w} \Big) < s < N .
\end{align*}
Then a function $f \in L_{\max(r,1)}^{loc} (\mathbb{R})$ belongs to $F^{s}_{r,q}(\mathbb{R})$ if and only if  
\begin{align*}
\|\, f \, |F^s_{r,q}(\mathbb{R})\|^{\clubsuit} := \|\, f \, |L_r (\mathbb{R})\| +
\Big\| \Big[
 \sum_{i= 0}^\infty 2^{iq(s + \frac{1}{w})} \Big(\int_{- 2^{-i}}^{2^{-i}} |\Delta_h^N f(x)|^{w} \, {dh} \Big)^{\frac{q}{w}} \Big]^{\frac{1}{q}}\Big|  L_r (\mathbb{R}) \Big\|
 <\infty .
\end{align*}
Furthermore $ \|\, f \, |F^s_{r,q}(\mathbb{R})\| $ and 
$ \|\, f \, |F^s_{r,q}(\mathbb{R})\|^{\clubsuit} $ are equivalent on $L_{\max(r,1)}^{loc} (\mathbb{R}).$
\end{lem}
Characterizations in terms of differences for the spaces $ F^s_{r,q}(\mathbb{R})   $ have a long history. We refer to Section 2.5.11 in \cite{Tr83}, Section 3.5.3 in \cite{Tr92} and Chapter 1.11.9 in \cite{Tr06}. The parameter $w$ in Lemma \ref{discrete} can be seen as an additional degree of freedom. However it influences the required conditions concerning the smoothness parameter $s$. 

\subsection{B-splines, quarks and quarklets}

In this section it is our main goal to define the quarklets. For that purpose in a first step we recall the definition of cardinal B-splines. The first order cardinal B-spline $  N_{1}  $ is just the characteristic function of the interval $ [0,1)  $, namely $ N_{1} := \chi_{[0,1)} $. Higher order cardinal B-splines of order $ m \in \mathbb{N}  $ with $ m \geq 2  $ are defined by induction using the convolution $ \ast $. So we have $ N_{m} := N_{m - 1} \ast N_{1}$. The cardinal B-splines possess some very nice properties. They are collected in the following lemma.

\begin{lem}\label{Bspl_elem}
Let $ m \in \mathbb{N}   $ and $ x \in \mathbb{R}  $. Then for the cardinal B-splines the following elementary properties hold. 

\begin{itemize}
\item[(i)] We can write $  N_{m}(x) = \frac{1}{(m-1)!} \sum_{k = 0}^{m}(-1)^{k} { m \choose k } (x - k)_{+}^{m-1}   $.

\item[(ii)] For $ m \geq 2  $ we have the recursion formula $ N_{m}(x) = \frac{x}{m-1} N_{m-1}(x) + \frac{m-x}{m-1} N_{m-1} (x-1)   $.

\item[(iii)] For $ m \geq 2  $ the derivatives satisfy $  N_{m}'(x) = N_{m-1}(x) - N_{m-1}(x-1)  $.

\item[(iv)] The B-splines are compactly supported with $ \supp N_{m} = [0,m]  $.

\item[(v)] We have $ \sum_{k \in \mathbb{Z} } N_{m}(x - k) = 1   $.

\end{itemize}
\end{lem}
Those properties are well-known. Let us refer to \cite{DeLo}, see Chapter 5.2 and 5.3.  One can also consult \cite{Ch1}. When we think on Triebel-Lizorkin spaces it is important to know under what conditions on the parameters the cardinal B-splines belong to the spaces $ F^{s}_{r,q}(\mathbb{R})   $. There is the following useful observation.

\begin{lem}\label{spline_inF}
Let $ s \in \mathbb{R}   $, $ 0 < r < \infty  $ and $ 0 < q \leq \infty  $. Let $ m \in \mathbb{N}  $. Then we have
\begin{align*}
N_{m} \in F^{s}_{r,q}(\mathbb{R}) \qquad \qquad \mbox{if and only if} \qquad \qquad s < m - 1 + \frac{1}{r} .
\end{align*}
\end{lem}
This result can be found in \cite{RS}, see Lemma 3 in Chapter 2.3.1. In what follows for fixed $ m \in \mathbb{N} $ we will work with the symmetrized cardinal B-spline $ \varphi(x)  :=   N_{m}  ( x + \lfloor \frac{m}{2} \rfloor   )    $. We observe $ \supp \varphi = [ - \lfloor \frac{m}{2} \rfloor    ,  \lceil \frac{m}{2} \rceil    ]   $. The symmetrized cardinal B-spline shows up in the following definition where we introduce the so-called quarks.

\begin{defi}\label{Bquark}
Let $ m \in \mathbb{N} $ and $ p \in \mathbb{N}_{0}  $. Then the p-th cardinal B-spline quark $ \varphi_{p}  $ is defined by
\begin{equation}
\varphi_{p}(x)  := \Big ( \frac{x}{\lceil \frac{m}{2} \rceil } \Big )^{p}  N_{m} \Big ( x + \lfloor \frac{m}{2} \rfloor  \Big ) .
\end{equation}
\end{defi}
The quarks will be very important for us in order to define the quarklets. Their properties have been studied in \cite{DaKRaa}. It is shown in \cite{CoDau} that for a given $ \tilde{m} \in \mathbb{N}  $ with $ \tilde{m} \geq m    $ and $  m + \tilde{m} \in 2 \mathbb{N}   $ there exists a compactly supported spline wavelet $ \psi  $ with
\begin{equation}\label{def_CDF_wav}
\psi = \sum_{k \in \mathbb{Z}} b_{k}  \varphi ( 2 \cdot - k   )
\end{equation}
with expansion coefficients $ b_{k} \in \mathbb{R}   $. Only finitely many of them are not zero. Moreover $ \psi  $ has $ \tilde{m}   $ vanishing moments and the  system
\begin{align*}
\Big \{ \varphi (  \cdot - k )  \ : \ k \in \mathbb{Z}  \Big  \} \cup \Big \{ 2^{\frac{j}{2}} \psi (2^{j} \cdot - k) \ : \ j \in \mathbb{N}_{0} \ , \ k \in \mathbb{Z}  \Big \}
\end{align*}
is a Riesz basis for $ L_{2}(\mathbb{R})   $. To construct such a $ \psi   $ we have to work with a compactly supported dual generator $ \tilde{\varphi}   $ associated to the primal generator $ \varphi $ that fulfills
\begin{equation}\label{biorto1}
\left \langle  \varphi , \tilde{\varphi} (\cdot - k)   \right\rangle_{L_{2}(\mathbb{R})} = \delta_{0,k} , \qquad k \in \mathbb{Z}.
\end{equation}
Connected with that there is another compactly supported biorthogonal wavelet $ \tilde{\psi} \in L_{2}(\mathbb{R})  $ with 
\begin{equation}\label{def_biort_wav1}
\tilde{\psi} = \sum_{k \in \mathbb{Z}}  \tilde{b}_{k} \tilde{\varphi} ( 2 \cdot - k   ).
\end{equation}
Here only finitely many of the $ \tilde{b}_{k} \in \mathbb{R}   $  are not zero. Moreover $ \tilde{\psi}  $ has $ m \in \mathbb{N}  $ vanishing moments and the system
\begin{align*}
\Big \{ \tilde{\varphi} (  \cdot - k )  \ : \ k \in \mathbb{Z}  \Big  \} \cup \Big \{ 2^{\frac{j}{2}} \tilde{\psi} (2^{j} \cdot - k) \ : \ j \in \mathbb{N}_{0} \ , \ k \in \mathbb{Z}  \Big \}
\end{align*}
is a Riesz basis for $ L_{2}(\mathbb{R})   $. For $ j \in \mathbb{N}_{0}   $ and $ k \in \mathbb{Z}   $ let us write
\begin{equation}\label{def_biort_wav2}
\psi_{j,k} = 2^{\frac{j}{2}} \psi ( 2^{j} \cdot - k ) \qquad \qquad \mbox{and} \qquad \qquad \tilde{\psi}_{j,k} = 2^{\frac{j}{2}} \tilde{\psi} ( 2^{j} \cdot - k ) .
\end{equation}
Moreover for $ k \in \mathbb{Z}   $ we put $ \psi_{-1,k} = \varphi (  \cdot - k )   $ and $ \tilde{\psi}_{-1,k} = \tilde{\varphi} (  \cdot - k )   $. Then we observe
\begin{equation}\label{biorto2}
\langle \psi_{j,k} , \tilde{\psi}_{j',k'}     \rangle_{L_{2}(\mathbb{R})} = \delta_{j , j'} \delta_{k , k'} , \qquad j, j' \in \mathbb{N}_{0}, \quad k, k' \in \mathbb{Z} .
\end{equation}
Moreover for each $ f \in L_{2}(\mathbb{R})  $ we have
\begin{align}
f & = \sum_{k \in \mathbb{Z}} \langle f , \tilde{\psi}_{-1,k}     \rangle_{L_{2}(\mathbb{R})} \psi_{-1,k}  + \sum_{j \in \mathbb{N}_{0} ,k \in \mathbb{Z}} \langle f , \tilde{\psi}_{j,k}     \rangle_{L_{2}(\mathbb{R})} \psi_{j,k} \nonumber \\
& = \sum_{k \in \mathbb{Z}} \langle f , \psi_{-1,k}     \rangle_{L_{2}(\mathbb{R})} \tilde{\psi}_{-1,k}  +  \sum_{j \in \mathbb{N}_{0}, k \in \mathbb{Z}} \langle f , \psi_{j,k}     \rangle_{L_{2}(\mathbb{R})} \tilde{\psi}_{j,k}    \label{def_biort_wav3}
\end{align}
with convergence in $ L_{2}(\mathbb{R})  $. For details and proofs concerning the above construction we refer to \cite{CoDau}, see especially Section 6.A. Now we can use our CDF-wavelets $ \psi $ to define the quarklets.  
\begin{defi}\label{def_quarklet}
Let $ p \in \mathbb{N}_{0}  $. Then the p-th quarklet $ \psi_{p} $ is defined by 
\begin{equation}
\psi_{p} := \sum_{k \in \mathbb{Z}} b_{k} \varphi_{p}(2 \cdot - k).
\end{equation}
Here the $ b_{k}  $ are the same as in \eqref{def_CDF_wav}. Moreover for $ j \in \mathbb{N}_{0}   $ and $ k \in \mathbb{Z}  $ we write 
\begin{equation}
\psi_{p,j,k} := 2^{\frac{j}{2}} \psi_{p}(2^{j} \cdot - k) \qquad \qquad \mbox{and} \qquad \qquad \psi_{p,-1,k} = \varphi_{p}( \cdot - k) .
\end{equation}
\end{defi}
Many properties of these quarklets as well as some first numerical applications already can be found in the literature. For that we refer to \cite{DaFKRaa}, \cite{DaKRaa} and \cite{DaRaaS}.

\section{Quarklet Characterizations for Triebel-Lizorkin spaces in the case $ s>0$}\label{sec_quas>0}

\subsection{Upper estimates}

As already mentioned the main goal of this paper is to prove quarklet characterizations for the Triebel-Lizorkin spaces $ F^s_{r,q}(\mathbb{R})  $. For that purpose in what follows we want to show some upper estimates for $ s > 0$. In a very first step we will deal with the special case $  r = q  $. Here we can use the fact that for the Besov spaces $ B^{s}_{r,r}(\mathbb{R})   $ already characterizations in terms of quarklets are known, see Section 4.5 and 4.6 in \cite{SiDiss}. Let $ \hat{p} \in \mathbb{N}_{0}    $ be fixed. Then we define the truncated quarklet system $ \Psi_{\hat{p}}  $ by 
\begin{equation}\label{trunc_qa_sys}
\Psi_{\hat{p}} := \{ \psi_{p,j,k} : 0 \leq p \leq \hat{p} , j \geq -1 , k \in \mathbb{Z}   \}.
\end{equation}
For these quarklet systems there is the following result concerning the Triebel-Lizorkin spaces $ F^s_{r,r}(\mathbb{R})   $.

\begin{prop}\label{thm_trunc_1}
Let $ m \in \mathbb{N}  $, $  \frac{1}{m} < r < \infty $  and $ \hat{p} \in \mathbb{N}_{0}   $. Let
\begin{equation}\label{scon_tr1B}
\max \Big(0, \frac 1r -1  \Big) < s < \min \Big ( m - 1 + \frac{1}{r} , m  \Big ).
\end{equation}
Let $ f \in \mathcal{S}'(\mathbb{R})  $ such that there exists a representation 
\begin{equation}\label{eq-inf1}
f = \sum_{p = 0}^{\hat{p}} \sum_{j \geq -1} \sum_{k \in \mathbb{Z}} c_{p,j,k} \psi_{p,j,k}
\end{equation}
with convergence in $ \mathcal{S}'(\mathbb{R})  $. Then there exists a $ C > 0  $ independent from $ f $ such that 
\begin{align*}
& \|\, f \, |F^s_{r,r}(\mathbb{R})\| \leq C \inf_{\eqref{eq-inf1}}  (\hat{p} + 1 )^{(2m+3)^{2}}  \Big \Vert \Big ( \sum_{p = 0}^{\hat{p}} \sum_{j \geq -1} \sum_{k \in \mathbb{Z}} 2^{jsr} 2^{j \frac{r}{2}} | c_{p,j,k}|^{r} | \chi_{j,k}( \cdot )   |^{r}  \Big )^{\frac{1}{r}}  \Big \vert L_{r}(\mathbb{R})   \Big \Vert .
\end{align*}
Here the infimum is taken over all sequences $ \{ c_{p,j,k}  \}_{p,j,k} \subset \mathbb{C}  $ such that \eqref{eq-inf1} is fulfilled.
\end{prop}

\begin{proof}
For the proof we use that we have $ F^s_{r,r}(\mathbb{R}) = B^s_{r,r}(\mathbb{R})    $ in the sense of equivalent quasi-norms, see for example Chapter 2.3.2. in \cite{Tr83}. Now take $ t \in \mathbb{R}  $ such that
\begin{equation}\label{scon_tr1B_t}
\max \Big(0, \frac 1r -1  \Big) < s < t < \min \Big ( m - 1 + \frac{1}{r} , m  \Big ).
\end{equation}
Because of \eqref{scon_tr1B} this is possible. Now we can apply the Theorems 4.26. and 4.31. from \cite{SiDiss} in slightly modified versions to obtain 
\begin{align*}
\|\, f \, |F^s_{r,r}(\mathbb{R})\| \leq C_{1} \inf_{\eqref{eq-inf1}} \max \Big( (\hat{p} + 1 )^{2t + 1 - \frac{1}{r}} , (\hat{p} + 1 )^{\frac{2tm}{r(m - 1 + \frac{1}{r})}} \Big) \Big ( \sum_{p,j,k}^{}   2^{j(s + \frac{1}{2} - \frac{1}{r})r}  | c_{p,j,k}|^{r}  \Big )^{\frac{1}{r}}  .
\end{align*}
Here the infimum is taken over all sequences $ \{ c_{p,j,k}  \}_{p,j,k}  $ such that \eqref{eq-inf1} is fulfilled. We use \eqref{scon_tr1B_t} to estimate
\begin{align*}
\max \Big( (\hat{p} + 1 )^{2t + 1 - \frac{1}{r}} , (\hat{p} + 1 )^{\frac{2tm}{r(m - 1 + \frac{1}{r})}} \Big) \leq  (\hat{p} + 1 )^{(2m+3)^{2}} .
\end{align*}
This estimate is not optimal but it is sufficient for our purposes. Later on it will be very important  that this estimate is
independent of the parameter $t$ since we want to use complex interpolation in another proof. Next we observe
\begin{align*}
& \Big \Vert \Big ( \sum_{p = 0}^{\hat{p}} \sum_{j \geq -1} \sum_{k \in \mathbb{Z}} 2^{jsr} 2^{j \frac{r}{2}} | c_{p,j,k} \chi_{j,k}( \cdot )   |^{r}  \Big )^{\frac{1}{r}}  \Big \vert L_{r}(\mathbb{R})   \Big \Vert \\
& \qquad \qquad = \Big (  \sum_{p = 0}^{\hat{p}} \sum_{j \geq -1} \sum_{k \in \mathbb{Z}} 2^{jsr} 2^{j \frac{r}{2}} | c_{p,j,k} |^{r}  \int_{- \infty}^{\infty} \chi_{j,k}( x )   dx \Big )^{\frac{1}{r}}.
\end{align*}
This in combination with $ \int_{- \infty}^{\infty} \chi_{j,k}( x )   dx = 2^{-j}  $ yields the desired result. 
\end{proof}
Next we turn to the full quarklet system $  \Psi $ given by
\begin{equation}
\Psi := \{ \psi_{p,j,k} : p \geq 0 , j \geq -1 , k \in \mathbb{Z}   \}.
\end{equation}
Here the counterpart of Proposition \ref{thm_trunc_1} reads as follows. 

\begin{prop}\label{stab_r=q_full}
Let $ m \in \mathbb{N}  $ and $  \frac{1}{m} < r < \infty $. Let 
\begin{equation}
\max \Big(0, \frac 1r -1  \Big) < s  < \min \Big ( m - 1 + \frac{1}{r} , m  \Big ).
\end{equation}
Let $ f \in \mathcal{S}'(\mathbb{R})  $ such that there exists a representation 
\begin{equation}\label{eq-inf2}
f = \sum_{p \geq 0} \sum_{j \geq -1} \sum_{k \in \mathbb{Z}} c_{p,j,k} \psi_{p,j,k}
\end{equation}
with convergence in $ \mathcal{S}'(\mathbb{R})  $. Then there exists a $ C > 0  $ independent from $ f $ such that 
\begin{align*}
 \|\, f \, |F^s_{r,r}(\mathbb{R})\| \leq C \inf_{\eqref{eq-inf2}} \Big \Vert \Big ( \sum_{p,j,k}  ( p + 1 )^{(2m+3)^{2} r}  2^{jsr} 2^{j \frac{r}{2}} | c_{p,j,k}|^{r}  | \chi_{j,k}   |^{r}  \Big )^{\frac{1}{r}}  \Big \vert L_{r}(\mathbb{R})   \Big \Vert .
\end{align*}
Here the infimum is taken over all sequences $ \{ c_{p,j,k}  \}_{p,j,k} \subset \mathbb{C} $ such that \eqref{eq-inf2} is fulfilled. 
\end{prop}
\begin{proof}
This result can be proved in the same way as Proposition \ref{thm_trunc_1}. We use the Theorems 4.27. and 4.32. from \cite{SiDiss}. Notice that again in a first step we have to work with an auxiliary parameter $ t $ given by \eqref{scon_tr1B_t}. With this we find
\begin{align*}
\|\, f \, |F^s_{r,r}(\mathbb{R})\|  \leq C \inf_{\eqref{eq-inf2}} \Big \Vert \Big ( \sum_{p,j,k}  \max \Big( ( p + 1 )^{(2t + 1 - \frac{1}{r} + \delta - \frac{\delta}{r})r} , ( p + 1 )^{\frac{2tmr}{r(m - 1 + \frac{1}{r})}} \Big) 2^{jsr} 2^{j \frac{r}{2}} | c_{p,j,k} \chi_{j,k}   |^{r}  \Big )^{\frac{1}{r}}  \Big \vert L_{r}(\mathbb{R})   \Big \Vert .
\end{align*}
Here $  \delta > 1 $ is another auxiliary parameter. In the proof of Theorem 4.27. in \cite{SiDiss} it is possible to choose $ \delta = 2  $. Then we observe
\begin{align*}
\max \Big( ( p + 1 )^{(2t + 1 - \frac{1}{r} + \delta - \frac{\delta}{r})r} , ( p + 1 )^{\frac{2tmr}{r(m - 1 + \frac{1}{r})}} \Big)  \leq  ( p + 1 )^{(2m+3)^{2} r}.
\end{align*}
This leads to the desired result.
\end{proof}
Now let us turn to the case that the parameters $r$ and $q$ are different. Then the situation becomes much more difficult. Therefore for a start we assume $ s  < m - 1   $. Then we find the following upper estimate.

\begin{prop}\label{res_Fspq_m-1}
Let $  0 < r < \infty $, $ 0 < q < \infty  $ and $ m \in \mathbb{N}  $ with $ m \geq 2  $. Let
\begin{equation}\label{res_Fspq_m-1_conds}
\max \Big(0, \frac 1r -1 , \frac 1q -1  \Big) < s  < m - 1.
\end{equation}
Let $ f \in \mathcal{S}'(\mathbb{R})  $ such that there exists a representation 
\begin{equation}\label{eq-inf3}
f = \sum_{p \geq 0} \sum_{j \geq -1} \sum_{k \in \mathbb{Z}} c_{p,j,k} \psi_{p,j,k}
\end{equation}
with convergence in $ \mathcal{S}'(\mathbb{R})  $. Then there exists a $ C > 0  $ independent from $ f $ such that 
\begin{align*}
& \|\, f \, |F^s_{r,q}(\mathbb{R})\| \leq C \inf_{\eqref{eq-inf3}} \Big\| \Big[ \sum_{p = 0}^{\infty}  \sum_{j= -1}^\infty \sum_{k \in \mathbb{Z}} ( p + 1 )^{(2m+3)^{2} q}  2^{jsq} 2^{\frac{jq}{2}}    |c_{p,j,k}|^{q}  |\chi_{j,k}(x)|^{q}   \Big]^{\frac{1}{q}}\Big|  L_{r} (\mathbb{R}) \Big\| .
\end{align*}
Here the infimum is taken over all sequences $ \{ c_{p,j,k}  \}_{p,j,k} \subset \mathbb{C}  $ such that \eqref{eq-inf3} is fulfilled.
\end{prop}

\begin{proof}
\textit{Step 1. Some preparations.} For the proof we deal with functions of the form
\begin{equation}\label{m-1_step1eq0}
f = \sum_{p \geq 0} f_{p} =  \sum_{p \geq 0} \sum_{j \geq -1} f_{p,j}  = \sum_{p \geq 0} \sum_{j \geq -1} \sum_{k \in \mathbb{Z}} c_{p,j,k} \psi_{p,j,k} .
\end{equation}
For such a function we want to prove the estimate
\begin{equation}\label{m-1_step1eq1}
\|\, f \, |F^s_{r,q}(\mathbb{R})\| \leq C \Big\| \Big[ \sum_{p = 0}^{\infty}  \sum_{j= -1}^\infty \sum_{k \in \mathbb{Z}} ( p + 1 )^{(2m+3)^{2} q}  2^{jsq} 2^{\frac{jq}{2}}    |c_{p,j,k}|^{q}  |\chi_{j,k}(x)|^{q}   \Big]^{\frac{1}{q}}\Big|  L_{r} (\mathbb{R}) \Big\|.
\end{equation}
Let $ t \in \mathbb{R} $ be near $s$ such that $ s < t < m - 1 $. Moreover let $ 0 < w < \min(1,r,q)  $. Since $ w < 1 $ and because of \eqref{res_Fspq_m-1_conds} we can use Lemma \ref{discrete} to get
\begin{equation}\label{m-1_step1eq2}
\|\, f \, |F^s_{r,q}(\mathbb{R})\| \leq C_{1} \|\, f \, |L_r (\mathbb{R})\| + C_{1}
\Big\| \Big[
 \sum_{i= 0}^\infty 2^{iq(s + \frac{1}{w})} \Big(\int_{- 2^{-i}}^{2^{-i}} |\Delta_h^N f(x)|^{w} \, {dh} \Big)^{\frac{q}{w}} \Big]^{\frac{1}{q}}\Big|  L_r (\mathbb{R}) \Big\| .
\end{equation}
Here we have $ N \in \mathbb{N}  $ with $ m - 1 \geq N > t > s  $ such that $N$ is as small as possible. In what follows we will investigate both terms in \eqref{m-1_step1eq2} separately. 

\textit{Step 2. An estimate for the Lebesgue-norm.}
Let us start with an investigation of the term $ \|\, f \, |L_r (\mathbb{R})\|  $. By using the definition of the function $ f $, see \eqref{m-1_step1eq0}, we observe
\begin{align*}
\|\, f \, |L_r (\mathbb{R})\| & = \Big \|\, \sum_{p \geq 0}^{} \sum_{j \geq -1} \sum_{k \in \mathbb{Z}} c_{p,j,k} \psi_{p,j,k} \, \Big |L_r (\mathbb{R}) \Big \| \\
& \leq \Big \|\, \sum_{p \geq 0}^{} \sum_{j \geq -1} \sum_{k \in \mathbb{Z}} |c_{p,j,k}|  |\psi_{p,j,k}| \, \Big |L_r (\mathbb{R}) \Big \| \\ 
& \leq C_{1} \Big \|\, \sum_{p \geq 0}^{} \sum_{j \geq -1} \sum_{k \in \mathbb{Z}} 2^{\frac{j}{2}} |c_{p,j,k}|  | \chi_{j,k}(\cdot) | \, \Big |L_r (\mathbb{R}) \Big \| .
\end{align*}
Here $ C_{1}  $ depends on $m$, but is independent from $f$. We used that the quarklets have a compact support, see Definition \ref{def_quarklet} and formula \ref{def_CDF_wav}. Moreover we used the boundedness of the quarklets, see in addition Definition \ref{Bquark}. Now let us consider the case $ q \leq 1  $. Here because of $ l_{q} \subset l_{1}  $ we observe  
\begin{align*}
\|\, f \, |L_r (\mathbb{R})\| & \leq C_{2} \Big \|\, \Big ( \sum_{p \geq 0}^{} \sum_{j \geq -1} \sum_{k \in \mathbb{Z}} 2^{\frac{jq}{2}} |c_{p,j,k}|^{q}  | \chi_{j,k}(\cdot) |^{q} \Big )^{\frac{1}{q}} \, \Big |L_r (\mathbb{R}) \Big \| \\
& \leq C_{2} \Big \|\, \Big ( \sum_{p \geq 0}^{} \sum_{j \geq -1} \sum_{k \in \mathbb{Z}} 2^{jsq} 2^{\frac{jq}{2}} |c_{p,j,k}|^{q}  | \chi_{j,k}(\cdot) |^{q} \Big )^{\frac{1}{q}} \, \Big |L_r (\mathbb{R}) \Big \| .
\end{align*}
Here we used $ s > 0  $. Next we deal with the case $ q > 1  $. Then we find $ q' \geq 1  $ such that $ \frac{1}{q} + \frac{1}{q'}=1   $. Let $ \delta > 1  $ be fixed. Then H\"olders inequality yields
\begin{align*}
& \|\, f \, |L_r (\mathbb{R})\| \\
&  \leq C_{2} \Big \|\, \sum_{p \geq 0}^{} \sum_{j \geq -1} \sum_{k \in \mathbb{Z}} (p+1)^{\delta} (p+1)^{- \delta} 2^{js} 2^{-js} 2^{\frac{j}{2}} |c_{p,j,k}|  | \chi_{j,k}(\cdot) | \, \Big |L_r (\mathbb{R}) \Big \| \\
& \leq C_{2} \Big \|\,  \Big ( \sum_{p \geq 0} \sum_{j \geq -1} (p+1)^{- \delta q'} 2^{-jsq'} \Big )^{\frac{1}{q'}} \Big ( \sum_{p \geq 0} \sum_{j \geq -1} (p+1)^{ \delta q} 2^{jsq}  2^{\frac{jq}{2}} \Big [   \sum_{k \in \mathbb{Z}}  |c_{p,j,k}|  | \chi_{j,k}(\cdot) | \Big ]^{q} \Big )^{\frac{1}{q}} \, \Big |L_r (\mathbb{R}) \Big \| \\
& \leq C_{3} \Big \|\,  \Big ( \sum_{p \geq 0} \sum_{j \geq -1} (p+1)^{ \delta q} 2^{jsq}  2^{\frac{jq}{2}} \Big [   \sum_{k \in \mathbb{Z}}  |c_{p,j,k}|  | \chi_{j,k}(\cdot) | \Big ]^{q} \Big )^{\frac{1}{q}} \, \Big |L_r (\mathbb{R}) \Big \| .
\end{align*}
We used $ s > 0  $ and $ \delta q' > 1  $. Another application of H\"olders inequality leads to 
\begin{align*}
\|\, f \, |L_r (\mathbb{R})\| & \leq C_{3} \Big \|\,  \Big ( \sum_{p \geq 0} \sum_{j \geq -1} (p+1)^{ \delta q} 2^{jsq}  2^{\frac{jq}{2}} \Big [   \sum_{k \in \mathbb{Z}}  |c_{p,j,k}|  | \chi_{j,k}(\cdot) |^{\frac{1}{q}} | \chi_{j,k}(\cdot) |^{\frac{1}{q'}} \Big ]^{q} \Big )^{\frac{1}{q}} \, \Big |L_r (\mathbb{R}) \Big \| \\
& \leq C_{3} \Big \|\,  \Big ( \sum_{j \geq -1}  \sum_{p \geq 0}  \sum_{k \in \mathbb{Z}} (p+1)^{ \delta q} 2^{jsq}  2^{\frac{jq}{2}}   |c_{p,j,k}|^{q}  | \chi_{j,k}(\cdot) |  \Big (   \sum_{k' \in \mathbb{Z}}     | \chi_{j,k'}(\cdot) | \Big )^{q - 1}  \Big )^{\frac{1}{q}} \, \Big |L_r (\mathbb{R}) \Big \| \\
& \leq C_{3}  \Big \|\,  \Big ( \sum_{j \geq -1}  \sum_{p \geq 0} \sum_{k \in \mathbb{Z}} (p+1)^{ \delta q}   2^{jsq}  2^{\frac{jq}{2}}   |c_{p,j,k}|^{q}  | \chi_{j,k}(\cdot) |  \Big )^{\frac{1}{q}} \, \Big |L_r (\mathbb{R}) \Big \| .
\end{align*}
Since $ \delta > 1  $ was arbitrary, this is what we want to have. 

\textit{Step 3. Estimates for differences.}
Now in \eqref{m-1_step1eq2} we deal with the term that contains differences. Since this step is very technical and comprehensive we split it into several substeps.
 
\textit{Substep 3.1. Discover truncated quarklet systems.}
At first we make some estimates with the goal to reach a situation where we can work with truncated quarklet systems, see \eqref{trunc_qa_sys}. Since $ w < 1 $ and $ w < q  $ we find
\begin{align*}
& \Big\| \Big[
 \sum_{i= 0}^\infty 2^{iq(s + \frac{1}{w})} \Big(\int_{- 2^{-i}}^{2^{-i}} |\Delta_h^N f(x)|^{w} \, {dh} \Big)^{\frac{q}{w}} \Big]^{\frac{1}{q}}\Big|  L_r (\mathbb{R}) \Big\| \\
& \qquad \leq \Big\| \Big[
 \sum_{i= 0}^\infty 2^{iq(s + \frac{1}{w})} \Big( \sum_{p \geq 0} \int_{- 2^{-i}}^{2^{-i}} |\Delta_h^N f_{p}(x)|^{w} \, {dh} \Big)^{\frac{q}{w}} \Big]^{\frac{1}{q}}\Big|  L_r (\mathbb{R}) \Big\|       \\
& \qquad \leq \Big\|  \Big \{  \sum_{p \geq 0} \Big[
 \sum_{i= 0}^\infty 2^{iq(s + \frac{1}{w})} \Big(  \int_{- 2^{-i}}^{2^{-i}} |\Delta_h^N f_{p}(x)|^{w} \, {dh} \Big)^{\frac{q}{w}} \Big]^{\frac{w}{q}} \Big \}^{\frac{1}{w}}  \Big|  L_r (\mathbb{R}) \Big\|   .
\end{align*}
Now let $ \delta > 1  $ be fixed. Since $ \frac{q}{w} > 1   $ there exists  $ b > 1  $ such that $ \frac{1}{\frac{q}{w}} + \frac{1}{b} = 1    $. Then we can use the H\"older inequality to get
\begin{align*}
& \Big\| \Big[
 \sum_{i= 0}^\infty 2^{iq(s + \frac{1}{w})} \Big(\int_{- 2^{-i}}^{2^{-i}} |\Delta_h^N f(x)|^{w} \, {dh} \Big)^{\frac{q}{w}} \Big]^{\frac{1}{q}}\Big|  L_r (\mathbb{R}) \Big\| \\
& \qquad \leq \Big\|  \Big \{  \sum_{p \geq 0} (p+1)^{- \frac{ \delta }{b}} (p+1)^{ \frac{\delta }{b}} \Big[
 \sum_{i= 0}^\infty 2^{iq(s + \frac{1}{w})} \Big(  \int_{- 2^{-i}}^{2^{-i}} |\Delta_h^N f_{p}(x)|^{w} \, {dh} \Big)^{\frac{q}{w}} \Big]^{\frac{w}{q}} \Big \}^{\frac{1}{w}}  \Big|  L_r (\mathbb{R}) \Big\|   \\
& \qquad \leq \Big\|  \Big \{ \Big ( \sum_{p \geq 0} (p+1)^{- \delta  } \Big )^{\frac{1}{b}} \Big (  \sum_{p \geq 0} (p+1)^{ \frac{\delta  q}{b w}} 
 \sum_{i= 0}^\infty 2^{iq(s + \frac{1}{w})} \Big(  \int_{- 2^{-i}}^{2^{-i}} |\Delta_h^N f_{p}(x)|^{w} \, {dh} \Big)^{\frac{q}{w}} \Big )^{\frac{w}{q}}  \Big \}^{\frac{1}{w}}  \Big|  L_r (\mathbb{R}) \Big\|   .
\end{align*}
Since $ \delta > 1   $ the first series is finite. Hence we obtain
\begin{align*}
& \Big\| \Big[
 \sum_{i= 0}^\infty 2^{iq(s + \frac{1}{w})} \Big(\int_{- 2^{-i}}^{2^{-i}} |\Delta_h^N f(x)|^{w} \, {dh} \Big)^{\frac{q}{w}} \Big]^{\frac{1}{q}}\Big|  L_r (\mathbb{R}) \Big\| \\
& \qquad \leq C_{4} \Big\|    \Big (  \sum_{p \geq 0} (p+1)^{ \frac{\delta (q-w)}{ w}} 
 \sum_{i= 0}^\infty 2^{iq(s + \frac{1}{w})} \Big(  \int_{- 2^{-i}}^{2^{-i}} |\Delta_h^N f_{p}(x)|^{w} \, {dh} \Big)^{\frac{q}{w}} \Big )^{\frac{1}{q}}    \Big|  L_r (\mathbb{R}) \Big\|   .
\end{align*}
Here we also used $ \frac{1}{b} = \frac{q - w}{q}   $. Now in what follows we have the possibility to work with the function $ f_{p}  $ when we deal with the differences. We will see that this has some advantages.

\textit{Substep 3.2. Work with differences.} Now we want to estimate the term that contains differences from above by some maximal functions. Since $ w \leq 1 $ we have $ |a+b|^{w} \leq |a|^{w} + |b|^{w}   $. Because of this with \eqref{m-1_step1eq0} we obtain
\begin{align*}
\Big(\int_{- 2^{-i}}^{2^{-i}} |\Delta_h^N f_{p}(x)|^{w} \, {dh} \Big)^{\frac{q}{w}} \leq  \Big( \sum_{j = -1}^{i-1} \int_{- 2^{-i}}^{2^{-i}} |\Delta_h^N f_{p,j}(x)|^{w} \, {dh} + \sum_{j = i}^{\infty} \int_{- 2^{-i}}^{2^{-i}} |\Delta_h^N f_{p,j}(x)|^{w} \, {dh} \Big)^{\frac{q}{w}}  .
\end{align*}
Take $ \varepsilon > 0  $ such that $ 0 < \varepsilon < \min(s,t-s)   $. Then we can write
\begin{align*}
& \Big(\int_{- 2^{-i}}^{2^{-i}} |\Delta_h^N f_{p}(x)|^{w} \, {dh} \Big)^{\frac{q}{w}} \\
& \qquad \leq \Big( \sum_{j = -1}^{i-1} \frac{2^{-(i-j)(-\varepsilon w)}}{2^{-(i-j)(-\varepsilon w)}} \int_{- 2^{-i}}^{2^{-i}} |\Delta_h^N f_{p,j}(x)|^{w} \, {dh} + \sum_{j = i}^{\infty} \frac{2^{(j-i)\varepsilon w}}{2^{(j-i)\varepsilon w}} \int_{- 2^{-i}}^{2^{-i}} |\Delta_h^N f_{p,j}(x)|^{w} \, {dh} \Big)^{\frac{q}{w}}  .
\end{align*}
We have
\begin{align*}
\sum_{j = -1}^{i-1} 2^{(i-j)(-\varepsilon w)} + \sum_{j = i}^{\infty} 2^{(j-i)(-\varepsilon w)} = \frac{2^{\varepsilon w}-2^{- \varepsilon w(i+1)}+1}{2^{\varepsilon w} - 1} = C_{\varepsilon } .
\end{align*}
Now for $  -1 \leq j \leq i - 1  $ we put $ \lambda_{j} =  C_{\varepsilon}^{-1} 2^{(i-j)(-\varepsilon w)}  $. For $ j \geq i   $ we write $ \lambda_{j} =  C_{\varepsilon}^{-1} 2^{(j-i)(-\varepsilon w)} $. With them we obtain a convex combination. Using our new notation we find
\begin{align*}
& \Big(\int_{- 2^{-i}}^{2^{-i}} |\Delta_h^N f_{p}(x)|^{w} \, {dh} \Big)^{\frac{q}{w}} \\
&  \leq C_{\varepsilon}^{\frac{q}{w}}  \Big( \sum_{j = -1}^{i-1} \lambda_{j} 2^{-(i-j)(-\varepsilon w)} \int_{- 2^{-i}}^{2^{-i}} |\Delta_h^N f_{p,j}|^{w} \, {dh} + \sum_{j = i}^{\infty} \lambda_{j} 2^{(j-i)\varepsilon w} \int_{- 2^{-i}}^{2^{-i}} |\Delta_h^N f_{p,j}|^{w} \, {dh} \Big)^{\frac{q}{w}}  .
\end{align*}
Now since $ w < q $ the function $ g(x) = x^{\frac{q}{w}}  $ is convex. Therefore we can apply the Jensen inequality to obtain
\begin{align*}
& \Big(\int_{- 2^{-i}}^{2^{-i}} |\Delta_h^N f_{p}(x)|^{w} \, {dh} \Big)^{\frac{q}{w}} \\
&  \leq C_{\varepsilon}^{\frac{q}{w}}  \Big( \sum_{j = -1}^{i-1} \lambda_{j} 2^{-(i-j)(-\varepsilon )q} \Big ( \int_{- 2^{-i}}^{2^{-i}} |\Delta_h^N f_{p,j}|^{w} \, {dh} \Big )^{\frac{q}{w}} + \sum_{j = i}^{\infty} \lambda_{j} 2^{(j-i)\varepsilon q} \Big ( \int_{- 2^{-i}}^{2^{-i}} |\Delta_h^N f_{p,j}|^{w} \, {dh} \Big )^{\frac{q}{w}} \Big)  .
\end{align*}
With $ \lambda_{j} \leq 1  $ we get
\begin{align*}
& \Big(\int_{- 2^{-i}}^{2^{-i}} |\Delta_h^N f_{p}(x)|^{w} \, {dh} \Big)^{\frac{q}{w}} \\
&  \leq C_{\varepsilon}^{\frac{q}{w}}  \Big( \sum_{j = -1}^{i-1}  2^{-(i-j)(-\varepsilon)q} \Big ( \int_{- 2^{-i}}^{2^{-i}} |\Delta_h^N f_{p,j}|^{w} \, {dh} \Big )^{\frac{q}{w}} + \sum_{j = i}^{\infty}  2^{(j-i)\varepsilon q} \Big ( \int_{- 2^{-i}}^{2^{-i}} |\Delta_h^N f_{p,j}|^{w} \, {dh} \Big )^{\frac{q}{w}} \Big)  .
\end{align*}
That leads to
\begin{align*}
& J := \sum_{i= 0}^\infty 2^{iq(s + \frac{1}{w})} \Big(\int_{- 2^{-i}}^{2^{-i}} |\Delta_h^N f_{p}(x)|^{w} \, {dh} \Big)^{\frac{q}{w}} \\
& \qquad \qquad  \leq C_{4}^{\frac{q}{w}} \sum_{j= -1}^\infty   \sum_{i = j}^{\infty} 2^{iq(s + \frac{1}{w})}  2^{-(i-j)(-\varepsilon) q} \Big ( \int_{- 2^{-i}}^{2^{-i}} |\Delta_h^N f_{p,j}(x)|^{w} \, {dh} \Big )^{\frac{q}{w}} \\
& \qquad \qquad \qquad \qquad + C_{4}^{\frac{q}{w}} \sum_{j= -1}^\infty \sum_{i=0}^{j-1} 2^{iq(s + \frac{1}{w})}  2^{(j-i)\varepsilon q} \Big ( \int_{- 2^{-i}}^{2^{-i}} |\Delta_h^N f_{p,j}(x)|^{w} \, {dh} \Big )^{\frac{q}{w}}   .
\end{align*}
Now we insert the additional parameter $t$. It is used to create an additional sum. The result is the following:
\begin{align*}
& J  \leq C_{4}^{\frac{q}{w}} \sum_{j= -1}^\infty   \sum_{i = j}^{\infty} 2^{iq(s + \frac{1}{w})}  2^{-(i-j)(-\varepsilon) q} 2^{-iq(t+\frac{1}{w})} \Big [ 2^{iq \frac{1}{w}} \sum_{b = i}^{\infty} 2^{bqt} \Big ( \int_{- 2^{-b}}^{2^{-b}} |\Delta_h^N f_{p,j}(x)|^{w} \, {dh} \Big )^{\frac{q}{w}} \Big ] \\
& \qquad \qquad \qquad \qquad + C_{4}^{\frac{q}{w}} \sum_{j= -1}^\infty \sum_{i=0}^{j-1} 2^{iq(s + \frac{1}{w})}  2^{(j-i)\varepsilon q} \Big ( \int_{- 2^{-i}}^{2^{-i}} |\Delta_h^N f_{p,j}(x)|^{w} \, {dh} \Big )^{\frac{q}{w}}   .
\end{align*}
Now at first we will look at the case $ i < j  $. Here we apply the well-known formula
\begin{equation}
\Delta^{N}_{h}f_{p,j}(x) = \sum_{l = 0}^{N} (-1)^{N-l} { N \choose l } f_{p,j} ( x + l h ).
\end{equation}
Using this we observe
\begin{align*}
& 
 2^{iq\frac{1}{w}} \Big(\int_{- 2^{-i}}^{2^{-i}} |\Delta_h^N f_{p,j}(x)|^{w} \, {dh} \Big)^{\frac{q}{w}}  \leq C_{5} \sum_{l = 0}^{N}   2^{iq\frac{1}{w}} \Big(\int_{- 2^{-i}}^{2^{-i}} | f_{p,j} ( x + l h ) |^{w} \, {dh} \Big)^{\frac{q}{w}}  .
\end{align*}
In the case of $ l = 0  $ we get
\begin{align*}
& 2^{iq\frac{1}{w}} \Big(\int_{- 2^{-i}}^{2^{-i}} | f_{p,j} ( x  ) |^{w} \, {dh} \Big)^{\frac{q}{w}}  \leq C_{6}  | f_{p,j} ( x  ) |^{q}  .
\end{align*}
For $ 1 \leq l \leq N  $ we find
\begin{align*}
&   2^{iq\frac{1}{w}} \Big(\int_{- 2^{-i}}^{2^{-i}} | f_{p,j} ( x + l h ) |^{w} \, {dh} \Big)^{\frac{q}{w}}  \leq   \Big( \frac{1}{2^{-i}}\int_{- 2^{-i}}^{2^{-i}} | f_{p,j} ( x + N h ) |^{w} \, {dh} \Big)^{\frac{q}{w}}.
\end{align*}
Let $ \textbf{M}   $ denote the Hardy-Littlewood-Maximal Operator. For a definition and some basic properties one may consult Chapter 1.2.3 in \cite{Tr83}. It is not difficult to see that we have
\begin{align*}
&   2^{iq\frac{1}{w}} \Big(\int_{- 2^{-i}}^{2^{-i}} | f_{p,j} ( x + l h ) |^{w} \, {dh} \Big)^{\frac{q}{w}}  \leq C_{7}  (( \textbf{M}| f_{p,j}|^{w})(x))^{\frac{q}{w}}.
\end{align*}
Recall that we have
\begin{align*}
| f_{p,j} ( x  ) |^{q} = ( | f_{p,j} ( x  ) |^{w} )^{\frac{q}{w}} \leq  (( \textbf{M}| f_{p,j}|^{w})(x))^{\frac{q}{w}}.
\end{align*}
So all in all for the case $ i < j    $ we obtain
\begin{equation}\label{eq-together1}
2^{iq\frac{1}{w}} \Big(\int_{- 2^{-i}}^{2^{-i}} |\Delta_h^N f_{p,j}(x)|^{w} \, {dh} \Big)^{\frac{q}{w}}  \leq C_{8} (( \textbf{M}| f_{p,j}|^{w})(x))^{\frac{q}{w}}.
\end{equation}
Next we turn to the case $ i \geq j + 1   $. Recall that for all $ x \in \mathbb{R}   $ and all $ j \in \mathbb{N}_{0} \cup \{ -1 \}   $ we can find an $ l \in \mathbb{Z}   $ such that we have $ x \in [2^{-j-1}l,2^{-j-1}(l+1))   $. Let us assume that we have $ x \in [2^{-j-1}l,2^{-j-1}(l+\frac{1}{2}))    $. Then the case $ x \in [2^{-j-1}(l+\frac{1}{2}),2^{-j-1}(l+1))    $ can be investigated with similar arguments. Now we can write
\begin{align*}
& 2^{iq \frac{1}{w}} \sum_{b= i}^{\infty}  2^{bqt}   \Big(\int_{- 2^{-b}}^{2^{-b}} |\Delta_h^N f_{p,j}(x)|^{w} \, {dh} \Big)^{\frac{q}{w}}  \\
& \qquad = 2^{iq \frac{1}{w}} \sum_{b= i}^{\infty}   2^{bqt }  \Big(\int_{\substack{|h| <  2^{-b}\\ |h| < \frac{x - 2^{-j-1}l}{2N} }} |\Delta_h^N f_{p,j}(x)|^{w} \, {dh}  + \int_{\substack{|h| <  2^{-b}\\ |h| > \frac{x - 2^{-j-1}l}{2N} }} |\Delta_h^N f_{p,j}(x)|^{w} \, {dh}  \Big)^{\frac{q}{w}} .
\end{align*}
For $ x \in [2^{-j-1}l,2^{-j-1}(l+\frac{1}{2})]   $ and $ |h| < \frac{x - 2^{-j-1}l}{2N}   $ because of the smoothness properties of the involved functions for some small $ \eta > 0  $ we observe 
\begin{align*}
|\Delta_h^N f_{p,j}(x)| \leq C_{9} |h|^{N} \max_{y \in [2^{-j-1}l + \eta,2^{-j-1}(l+1) - \eta] } | f^{(N)}_{p,j}(y)   |.
\end{align*}
Now we benefit from the special structure of the functions $ f_{p,j}  $, see \eqref{m-1_step1eq0}. Recall that the quarklets consist of quarks that are built out of B-splines multiplied with polynomials, see Definition \ref{Bquark}. Therefore we can apply the Markov inequality (see for example Chapter 4.9.6 in \cite{Tim}) and get 
\begin{equation}\label{markov_inpr}
|\Delta_h^N f_{p,j}(x)| \leq C_{10} (p + 1)^{2N} 2^{jN} |h|^{N} \max_{y \in [2^{-j-1}l,2^{-j-1}(l+1)] } | f_{p,j}(y)   |.
\end{equation}
This leads to 
\begin{align*}
& \Big(\int_{\substack{|h| <  2^{-b}\\ |h| < \frac{x - 2^{-j-1}l}{2N} }} |\Delta_h^N f_{p,j}(x)|^{w} \, {dh} \Big)^{\frac{q}{w}} \\
& \qquad \qquad \leq C_{11} \max_{y \in [2^{-j-1}l,2^{-j-1}(l+1)] } | f_{p,j}(y)   |^{q} (p + 1)^{2Nq} 2^{jNq} \Big(\int_{|h| <  2^{-b}}  |h|^{Nw}  \, {dh} \Big)^{\frac{q}{w}} \\
& \qquad \qquad \leq C_{12} \max_{y \in [2^{-j-1}l,2^{-j-1}(l+1)] } | f_{p,j}(y)   |^{q} (p + 1)^{2Nq} 2^{jNq} 2^{-bNq} 2^{-bq\frac{1}{w}} .
\end{align*}
Consequently we also find
\begin{align*}
& 2^{iq \frac{1}{w}}  \sum_{b= i}^{\infty}   2^{bqt }  \Big(\int_{\substack{|h| <  2^{-b}\\ |h| < \frac{x - 2^{-j-1}l}{2N} }} |\Delta_h^N f_{p,j}(x)|^{w} \, {dh} \Big)^{\frac{q}{w}} \\
& \qquad \qquad \leq C_{13} 2^{iq \frac{1}{w}} \sum_{b= i}^{\infty}   2^{bqt }   \max_{y \in [2^{-j-1}l,2^{-j-1}(l+1)] } | f_{p,j}(y)   |^{q} (p + 1)^{2Nq} 2^{jNq} 2^{-bNq} 2^{-bq\frac{1}{w}} \\
& \qquad \qquad \leq C_{14}   \max_{y \in [2^{-j-1}l,2^{-j-1}(l+1)] } | f_{p,j}(y)   |^{q} (p + 1)^{2Nq} 2^{jNq} 2^{iq \frac{1}{w}} \sum_{b= i}^{\infty}   2^{bq(t-N) } 2^{-bq\frac{1}{w}}  .
\end{align*}
For $ b \geq i $ we have $ 2^{-bq\frac{1}{w}} \leq 2^{-iq\frac{1}{w}}   $. Moreover since $ t < N  $ for the geometric series we calculate
\begin{align*}
\sum_{b= i}^{\infty}  2^{bq(t-N) } \leq C_{15} \frac{2^{iq(t-N)}}{1 - 2^{q(t-N)}} \leq C_{16} 2^{iq(t-N)}.
\end{align*}
Then with $ i \geq j  $ we get $ 2^{iq(t-N)} \leq 2^{jq(t-N)}   $. Therefore we have 
\begin{align*}
& 2^{iq \frac{1}{w}} \sum_{b= i}^{\infty}   2^{bqt }  \Big(\int_{\substack{|h| <  2^{-b}\\ |h| < \frac{x - 2^{-j-1}l}{2N} }} |\Delta_h^N f_{p,j}(x)|^{w} \, {dh} \Big)^{\frac{q}{w}} \\
& \qquad \qquad \leq C_{17}   \max_{y \in [2^{-j-1}l,2^{-j-1}(l+1)] } | f_{p,j}(y)   |^{q} (p + 1)^{2Nq}  2^{jqt}    \\
& \qquad \qquad \leq C_{18}   (p + 1)^{(N+\frac{1}{w})2q}  2^{jqt}   (( \textbf{M}| f_{p,j}|^{w})(x))^{\frac{q}{w}}    .
\end{align*}
In the last step to switch from the maximum to the maximal function we used the special structure of $ f_{p,j}  $ in combination with an estimate from page 236 in \cite{Tim}. Now we turn to the case $   |h| > \frac{x - 2^{-j-1}l}{2N}   $ and $  |h| <  2^{-b}  $. Because of $ j \leq i \leq b  $ we have $  |h| <  2^{-b} \leq 2^{-j}   $. Again let $ x \in [2^{-j-1}l,2^{-j-1}(l+\frac{1}{2})]   $. For what follows it is very important that we have $ s < t < N \leq m - 1  $. Then because of the smoothness properties of the quarklets we observe
\begin{equation}\label{pro_unsm_1}
|\Delta_h^N f_{p,j}(x)|  \leq  C_{19} |h|^{N-1} \max_{y \in [2^{-j}(\frac{1}{2}l -(N-1)),2^{-j}(\frac{1}{2}(l+\frac{1}{2})+(N-1))] } |\Delta_h^1 f^{(N-1)}_{p,j}(y)   | .
\end{equation}
Notice that the quarklets are not smooth enough to work with $ f^{(N)}_{p,j}  $ on the whole interval $ [2^{-j}(\frac{1}{2}l - N),2^{-j}(\frac{1}{2}(l+\frac{1}{2})+ N)]    $. Therefore this time we have to go another way. We have  
\begin{align*}
|\Delta_h^N f_{p,j}(x)| & \leq  C_{20} |h|^{N-1} \max_{y \in [2^{-j}(\frac{1}{2}l -(N-1)),2^{-j}(\frac{1}{2}(l+\frac{1}{2})+(N-1))] } | f^{(N-1)}_{p,j}(y + h) - f^{(N-1)}_{p,j}(y)  | \\
& =  C_{20}  |h|^{N-1}  | f^{(N-1)}_{p,j}(y^{*} + h) - f^{(N-1)}_{p,j}(y^{*})  | .
\end{align*}
Here $ y^{*} \in  [2^{-j}(\frac{1}{2}l -(N-1)),2^{-j}(\frac{1}{2}(l+\frac{1}{2})+(N-1))]  $ is such that the maximum is realized. Since we are working with continuous functions on compact intervals such a number exists. Now let $ y_{crit}  $ be a critical point where the function is not differentiable between $ y^{*} $ and $ y^{*} + h  $. Because of $  |h| <  2^{-b}  $ there can be only one of these points. We can write
\begin{align*}
& | f^{(N-1)}_{p,j}(y^{*} + h) - f^{(N-1)}_{p,j}(y^{*})  | \\
& \qquad \leq | f^{(N-1)}_{p,j}(y^{*} + h) - f^{(N-1)}_{p,j}(y_{crit}) | + | f^{(N-1)}_{p,j}(y_{crit}) - f^{(N-1)}_{p,j}(y^{*})  | .
\end{align*}
So we found two intervals where we can work with smooth functions. The Mean Value Theorem yields 
\begin{align*}
& | f^{(N-1)}_{p,j}(y^{*} + h) - f^{(N-1)}_{p,j}(y^{*})  |  \leq |h| | f^{(N)}_{p,j}(\varrho_{1}) | + |h| | f^{(N)}_{p,j}(\varrho_{2})  | .
\end{align*}
Here $ \varrho_{1} \in \mathbb{R}  $ is between $ y^{*} + h   $ and $   y_{crit}$. The number $ \varrho_{2} \in \mathbb{R}  $ is between $ y^{*}    $ and $   y_{crit}$. When we combine this with \eqref{pro_unsm_1} and use the Markov inequality as in formula \eqref{markov_inpr} we get
\begin{align*}
|\Delta_h^N f_{p,j}(x)| & \leq  C_{21} (p + 1)^{2N } 2^{jN} |h|^{N} \max_{y \in [2^{-j}(\frac{1}{2}l - N),2^{-j}(\frac{1}{2}(l+\frac{1}{2})+ N)] } | f_{p,j}(y) |.
\end{align*}
Notice that this is similar to \eqref{markov_inpr}. Therefore we can proceed like there to obtain
\begin{equation}\label{eq-together2}
2^{iq \frac{1}{w}} \sum_{b= i}^{\infty}   2^{bqt }  \Big(\int_{\substack{|h| <  2^{-b}\\ |h| > \frac{x - 2^{-j-1}l}{2N} }} |\Delta_h^N f_{p,j}(x)|^{w} \, {dh} \Big)^{\frac{q}{w}}  \leq C_{22}   (p + 1)^{(N+\frac{1}{w})2q}  2^{jqt}   (( \textbf{M}| f_{p,j}|^{w})(x))^{\frac{q}{w}}    .
\end{equation}
Now we put everything together. So a combination of \eqref{eq-together1} and \eqref{eq-together2} yields
\begin{align*}
& J  \leq C_{23} \sum_{j= -1}^\infty   \sum_{i = j}^{\infty} 2^{iqs}  2^{-(i-j)(-\varepsilon) q} 2^{-iqt} \Big [  (p + 1)^{(N+\frac{1}{w})2q}  2^{jqt}   (( \textbf{M}| f_{p,j}|^{w})(x))^{\frac{q}{w}}   \Big ] \\
& \qquad \qquad  + C_{23} \sum_{j= -1}^\infty \sum_{i=0}^{j-1} 2^{iqs }  2^{(j-i)\varepsilon q}   (( \textbf{M}| f_{p,j}|^{w})(x))^{\frac{q}{w}}. 
\end{align*}
Due to $  \varepsilon < \min(s,t-s)   $ the series converge. Hence we find
\begin{align*}
 J &  \leq C_{24} \sum_{j= -1}^\infty 2^{jsq}   (p + 1)^{(N+\frac{1}{w})2q}    (( \textbf{M}| f_{p,j}|^{w})(x))^{\frac{q}{w}}   \sum_{i = j}^{\infty}    2^{-(i-j)(t-\varepsilon-s)q} \\
&  \qquad \qquad + C_{24} \sum_{j= -1}^\infty 2^{jsq}   (( \textbf{M}| f_{p,j}|^{w})(x))^{\frac{q}{w}} \sum_{i=0}^{j-1} 2^{(j-i)(\varepsilon - s) q}  \\
&   \leq C_{25} \sum_{j= -1}^\infty 2^{jsq}     (p + 1)^{(N+\frac{1}{w})2q}   (( \textbf{M}| f_{p,j}|^{w})(x))^{\frac{q}{w}} . 
\end{align*}
So we were able to find the desired estimate in terms of maximal functions.

\textit{Substep 3.3. Use the special form of $f$ to complete the proof.} Now we are able to complete the proof. For that purpose at first we have to combine the results from substep 3.1 and substep 3.2. We observe
\begin{align*}
& \Big\| \Big[
 \sum_{i= 0}^\infty 2^{iq(s + \frac{1}{w})} \Big(\int_{- 2^{-i}}^{2^{-i}} |\Delta_h^N f(x)|^{w} \, {dh} \Big)^{\frac{q}{w}} \Big]^{\frac{1}{q}}\Big|  L_r (\mathbb{R}) \Big\| \\
& \qquad \leq C_{26} \Big\|    \Big (  \sum_{p \geq 0} \sum_{j= -1}^\infty (p+1)^{ \frac{\delta (q-w)}{ w} + (N+\frac{1}{w})2q}   2^{jsq}       (( \textbf{M}| f_{p,j}|^{w})(x))^{\frac{q}{w}} 
  \Big )^{\frac{1}{q}}    \Big|  L_r (\mathbb{R}) \Big\| \\
& \qquad = C_{26} \Big\| \Big[ \sum_{p \geq 0}  \sum_{j= -1}^\infty  (( \textbf{M}| 2^{js}   (p+1)^{ \frac{\delta (q-w)}{ wq} + (N+\frac{1}{w})2}    f_{p,j}|^{w})(x))^{\frac{q}{w}}  \Big]^{\frac{w}{q}}\Big|  L_{\frac{r}{w}} (\mathbb{R}) \Big\|^{\frac{1}{w}} .
\end{align*}
Recall that $  0 < w < \min(1,r,q)  $. Therefore we have $ \frac{q}{w} > 1   $ and $ \frac{r}{w} > 1   $. Hence we can apply the Hardy-Littlewood-Maximal Inequality. For that we refer to Chapter 1.2.3 in \cite{Tr83}. We find
\begin{align*}
& \Big\| \Big[
 \sum_{i= 0}^\infty 2^{iq(s + \frac{1}{w})} \Big(\int_{- 2^{-i}}^{2^{-i}} |\Delta_h^N f(x)|^{w} \, {dh} \Big)^{\frac{q}{w}} \Big]^{\frac{1}{q}}\Big|  L_r (\mathbb{R}) \Big\| \\
& \qquad \leq C_{27} \Big\| \Big[ \sum_{p \geq 0}  \sum_{j= -1}^\infty 2^{jsq}    (p+1)^{ \frac{\delta (q-w)}{ w} + (N+\frac{1}{w})2q}      | f_{p,j}(x)|^{q}  \Big]^{\frac{1}{q}}\Big|  L_{r} (\mathbb{R}) \Big\| .
\end{align*}
Now we plug in the definition of the functions $ f_{p,j}  $. When we use the definition of the quarklets, see Definition \ref{def_quarklet}, we obtain 
\begin{align*}
& \Big\| \Big[
 \sum_{i= 0}^\infty 2^{iq(s + \frac{1}{w})} \Big(\int_{- 2^{-i}}^{2^{-i}} |\Delta_h^N f(x)|^{w} \, {dh} \Big)^{\frac{q}{w}} \Big]^{\frac{1}{q}}\Big|  L_r (\mathbb{R}) \Big\| \\
& \qquad \leq C_{28} \Big\| \Big[ \sum_{p \geq 0} \sum_{j= -1}^\infty 2^{jsq}    (p+1)^{ \frac{\delta (q-w)}{ w} + (N+\frac{1}{w})2q}   \Big |   \sum_{k \in \mathbb{Z}} c_{p,j,k} \psi_{p,j,k}(x) \Big |^{q}  \Big]^{\frac{1}{q}}\Big|  L_{r} (\mathbb{R}) \Big\| \\
& \qquad \leq C_{29} \Big\| \Big[ \sum_{p \geq 0}  \sum_{j= -1}^\infty 2^{jsq} 2^{\frac{jq}{2}}   (p+1)^{ \frac{\delta (q-w)}{ w} + (N+\frac{1}{w})2q}   \Big (   \sum_{k \in \mathbb{Z}} |c_{p,j,k}|  |\chi_{j,k}(x)| \Big )^{q}  \Big]^{\frac{1}{q}}\Big|  L_{r} (\mathbb{R}) \Big\| .
\end{align*}
In the last step we used the support properties of the quarklets. Next we want to move the sum concerning $ k $. In the case $ q \leq 1  $ because of $ | a + b |^{q} \leq |a|^{q} + |b|^{q}      $ this is simple. For $ q > 1  $ we find a $ q'  $ such that $ \frac{1}{q} + \frac{1}{q'} = 1    $. We write $ |\chi_{j,k}(x)| = |\chi_{j,k}(x)|^{\frac{1}{q}}  |\chi_{j,k}(x)|^{\frac{1}{q'}}    $. Then with the H\"older inequality we get
\begin{align*}
& \Big\| \Big[
 \sum_{i= 0}^\infty 2^{iq(s + \frac{1}{w})} \Big(\int_{- 2^{-i}}^{2^{-i}} |\Delta_h^N f(x)|^{w} \, {dh} \Big)^{\frac{q}{w}} \Big]^{\frac{1}{q}}\Big|  L_r (\mathbb{R}) \Big\| \\
& \qquad \leq C_{30} \Big\| \Big[ \sum_{p \geq 0}  \sum_{j= -1}^\infty \sum_{k \in \mathbb{Z}}  2^{jsq} 2^{\frac{jq}{2}}   (p+1)^{ \frac{\delta (q-w)}{ w} + (N+\frac{1}{w})2q}    |c_{p,j,k}|^{q}  |\chi_{j,k}(x)|^{q}   \Big]^{\frac{1}{q}}\Big|  L_{r} (\mathbb{R}) \Big\| .
\end{align*}
Now recall that we have $ 0 < w < \min(1,r,q)  $. For the proof it is not necessary to choose $w$ near zero. It is enough if $w$ is only a little bit smaller than $ \min(1,r,q)  $. For example $ w = \frac{1}{2}  \min(1,r,q)  $ works. Moreover we used $ \delta > 1  $. Therefore it is possible to put $ \delta = 2  $. When we use this in combination with \eqref{res_Fspq_m-1_conds} we obtain
\begin{align*}
(p+1)^{ \frac{\delta (q-w)}{ w} + (N+\frac{1}{w})2q} \leq ( p + 1 )^{[ 2m \delta + 6m + 1]q} \leq ( p + 1 )^{(2m+3)^{2} q}.
\end{align*}
Of course this estimate is not optimal. But since in what follows we want to work uniformly, in view of Proposition \ref{stab_r=q_full} we can not be much better here. This also leads to
\begin{align*}
& \Big\| \Big[
 \sum_{i= 0}^\infty 2^{iq(s + \frac{1}{w})} \Big(\int_{- 2^{-i}}^{2^{-i}} |\Delta_h^N f(x)|^{w} \, {dh} \Big)^{\frac{q}{w}} \Big]^{\frac{1}{q}}\Big|  L_r (\mathbb{R}) \Big\| \\
& \qquad \leq C_{31}  \Big\| \Big[ \sum_{p = 0}^{\infty}  \sum_{j= -1}^\infty \sum_{k \in \mathbb{Z}} ( p + 1 )^{(2m+3)^{2} q}  2^{jsq} 2^{\frac{jq}{2}}    |c_{p,j,k}|^{q}  |\chi_{j,k}(x)|^{q}   \Big]^{\frac{1}{q}}\Big|  L_{r} (\mathbb{R}) \Big\| .
\end{align*}
Since the representation was arbitrary, taking the infimum gives the claim. The proof is complete. 
\end{proof}

Our next goal is to prove a counterpart of Proposition \ref{res_Fspq_m-1} that also holds for $ s \geq m - 1  $. For that purpose we want to use complex interpolation. Let $(X_0,X_1)$ be an interpolation couple of quasi-Banach spaces. Then by $[X_0,X_1]_\theta$ we denote the result of the complex interpolation of these spaces. We refer to Calder{\'o}n \cite{ca64}, Bergh, L\"ofstr\"om \cite{BL},
Kre{\u{\i}}n, Petunin, Semenov  \cite{KPS}, Lunardi \cite{lun} and Triebel
\cite{t78} for a detailed discussion. For the Triebel-Lizorkin spaces the behavior with respect to complex interpolation is already known. We recall the following result. 

\begin{lem}\label{lem_comint_Fsrq1}
Let $ s_{0}, s_{1} \in \mathbb{R}  $, $ 0 < r_{0}, r_{1} < \infty $ and $  0 < q_{0},q_{1} < \infty $. Let $ 0 < \theta < 1 $ and put 
\begin{equation}\label{inter-seqor1}
s = ( 1 - \theta ) s_{0} + \theta s_{1} , \qquad  \qquad \frac{1}{r} = \frac{1 - \theta}{r_{0}} + \frac{\theta}{r_{1}} , \qquad  \qquad \frac{1}{q} = \frac{1 - \theta}{q_{0}} + \frac{\theta}{q_{1}} .
\end{equation}
Then we have
\begin{align*}
[ F^{s_{0}}_{r_{0},q_{0}}(\mathbb{R})  , F^{s_{1}}_{r_{1},q_{1}}(\mathbb{R}) ]_\theta  =  F^s_{r,q}(\mathbb{R})  .
\end{align*}
\end{lem}
For the case of quasi-Banach spaces this result can be found in \cite{KMM}, see Theorem 9.1. For the case of Banach spaces this result is known since a long time. Here we refer to Section 2.4. in \cite{t78} or to Corollary 8.3 in \cite{FJ1}. Notice that the complex interpolation method $[X_0,X_1]_\theta $ has the so-called interpolation property, and this also holds for quasi-Banach spaces. For explanations and proofs concerning that topic we refer to \cite{KMM}, see especially Theorem 8.1. To prove a counterpart of Proposition \ref{res_Fspq_m-1} we need some more knowledge concerning complex interpolation. More precisely we have to investigate the following modified sequence spaces.  

\begin{defi}\label{def_seqsp_qua}
Let $ 0 < r < \infty    $, $ 0 < q \leq \infty  $ and $ s \in \mathbb{R}  $ with $ s \not = 0  $. Moreover let $ m \in \mathbb{N}  $. By $ \lambda = \{ c_{p,j,k} \in \mathbb{C} : p \in \mathbb{N}_{0}, j \in \mathbb{N}_{0} \cup \{-1 \} , k \in \mathbb{Z}   \} $ we denote a complex sequence. Then the sequence space $ f^{s}_{r,q}(m)   $ is the collection of all sequences $  \lambda  $ such that
\begin{align*}
\Vert \lambda \vert f^{s}_{r,q}(m)   \Vert :=   \Big\| \Big[ \sum_{p = 0}^{\infty}  \sum_{j= -1}^\infty \sum_{k \in \mathbb{Z}}  ( p + 1 )^{ \sgn(s) (2m+3)^{2} q}  2^{jsq} 2^{\frac{jq}{2}}    |c_{p,j,k}|^{q}  |\chi_{j,k}(x)|^{q}   \Big]^{\frac{1}{q}}\Big|  L_{r} (\mathbb{R}) \Big\| < \infty . 
\end{align*}
\end{defi}
We want to prove a counterpart of Lemma \ref{lem_comint_Fsrq1} for our sequence spaces $ f^{s}_{r,q}(m)   $. Let us remark that for the original sequence spaces without polynomial enrichment such results are already known, see Theorem 9.1 and formula (9.12) in \cite{KMM}. As a preparation we have to deal with Calderon products. They are defined in the following way.

\begin{defi}\label{def_cald_P}
Let $ (X_{0} , \Vert \cdot \vert X_{0} \Vert )    $ and $ (X_{1} , \Vert \cdot \vert X_{1} \Vert )    $ be two quasi-Banach lattices. Let $ 0 < \theta < 1$. Then the Calderon product $ X_{0}^{1 - \theta}  X_{1}^{ \theta} $ is the collection of all measurable functions $h$ such that there exist $  f \in X_{0} $ with $ \Vert f \vert X_{0} \Vert \leq 1 $ and $  g \in X_{1}  $ with $ \Vert g \vert X_{1} \Vert \leq 1 $ such that for a $ \lambda > 0  $ we have 
\begin{equation}\label{cald_p_lambda}
 |h| \leq \lambda |f|^{1 - \theta} |g|^{\theta}   .
\end{equation}
Moreover we define the quasi-norm
\begin{align*}
\Vert h \vert  X_{0}^{1 - \theta}  X_{1}^{ \theta} \Vert := \inf \{ \lambda > 0 \ : \ \eqref{cald_p_lambda} \ \mbox{holds with} \ \Vert f \vert X_{0} \Vert \leq 1  \ \mbox{and} \ \Vert g \vert X_{1} \Vert \leq 1  \} .
\end{align*}
\end{defi}

In what follows we want to compute the Calderon product for two of our sequence spaces $  f^{s}_{r,q}(m)     $. There is the following result.

\begin{lem}\label{lem_cald_seq}
Let $ s_{0}, s_{1} > 0  $, $ 0 < r_{0}, r_{1} < \infty $ and $  0 < q_{0},q_{1} < \infty $. Let $ m \in \mathbb{N}  $. Let $ 0 < \theta < 1 $ and put 
\begin{align*}
s = ( 1 - \theta ) s_{0} + \theta s_{1} , \qquad  \qquad \frac{1}{r} = \frac{1 - \theta}{r_{0}} + \frac{\theta}{r_{1}} , \qquad  \qquad \frac{1}{q} = \frac{1 - \theta}{q_{0}} + \frac{\theta}{q_{1}} .
\end{align*}
Then we have $  f^{s}_{r,q}(m)  =  ( f^{s_{0}}_{r_{0},q_{0}}(m) )^{1 - \theta}    ( f^{s_{1}}_{r_{1},q_{1}}(m) )^{ \theta} $ in the sense of equivalent quasi-norms. 
\end{lem} 

\begin{proof}
For the proof we use some ideas from Theorem 8.2 from \cite{FJ2}.

\textit{Step 1.} At first let $  \lambda = \{ c_{p,j,k}   \}_{p,j,k}   $ be a sequence with $  \lambda \in ( f^{s_{0}}_{r_{0},q_{0}}(m) )^{1 - \theta}    ( f^{s_{1}}_{r_{1},q_{1}}(m) )^{ \theta}    $. Let $ \varepsilon > 0  $ and put $ B = (1 + \varepsilon) \Vert \lambda \vert  ( f^{s_{0}}_{r_{0},q_{0}}(m) )^{1 - \theta}    ( f^{s_{1}}_{r_{1},q_{1}}(m) )^{ \theta}   \Vert    $. Then there exist squences $ \alpha = \{ \alpha_{p,j,k}   \}_{p,j,k} $ and $ \beta = \{ \beta_{p,j,k}   \}_{p,j,k} $ with $ \Vert \alpha \vert  f^{s_{0}}_{r_{0},q_{0}}(m)   \Vert \leq 1    $ and $ \Vert \beta \vert  f^{s_{1}}_{r_{1},q_{1}}(m)   \Vert \leq 1    $ such that we have
\begin{align*}
|  c_{p,j,k}   | \leq B |  \alpha_{p,j,k} |^{1 - \theta} |  \beta_{p,j,k} |^{\theta}
\end{align*}
for all $ p \in \mathbb{N}_{0} $, $ j \geq -1 $ and $ k \in \mathbb{Z}   $. This follows from Definition \ref{def_cald_P}. Next we observe $ \frac{(1 - \theta )q}{q_{0}} + \frac{q \theta}{q_{1}} = 1 $. Therefore we can use the H\"older inequality to get
\begin{align*}
\Vert \lambda \vert f^{s}_{r,q}(m)   \Vert & \leq B  \Big\|  \Big ( \sum_{p,j,k}   ( p + 1 )^{  (2m+3)^{2} q_{0}}  2^{\frac{j q_{0}}{2}} 2^{j q_{0} s_{0}  }      |  \alpha_{p,j,k} |^{ q_{0}  }    |\chi_{j,k}(x)|  \Big )^{\frac{(1 - \theta )}{q_{0}}}        \\
& \qquad \qquad \qquad  \times \Big (  \sum_{p,j,k}    ( p + 1 )^{  (2m+3)^{2} q_{1} }       2^{\frac{j q_{1} }{2}} 2^{j s_{1} q_{1} }   |  \beta_{p,j,k} |^{ q_{1} }  |\chi_{j,k}(x)|  \Big )^{\frac{ \theta}{q_{1}}}      \Big|  L_{r} (\mathbb{R}) \Big\|  . 
\end{align*}
For what follows we use $  \frac{(1 - \theta )r}{r_{0}} + \frac{r \theta}{r_{1}} = 1     $. Then the H\"older inequality yields
\begin{align*}
\Vert \lambda \vert f^{s}_{r,q}(m)   \Vert & \leq B   \Big ( \int_{\mathbb{R}}  \Big ( \sum_{p,j,k}   ( p + 1 )^{  (2m+3)^{2} q_{0}}  2^{\frac{j q_{0}}{2}} 2^{j q_{0} s_{0}  }      |  \alpha_{p,j,k} |^{ q_{0}  }    |\chi_{j,k}(x)|  \Big )^{\frac{r_{0}}{q_{0}}} dx \Big )^{\frac{(1 - \theta )}{r_{0}}}       \\
&  \qquad \qquad \qquad \times  \Big ( \int_{\mathbb{R}}  \Big (  \sum_{p,j,k}    ( p + 1 )^{  (2m+3)^{2} q_{1} }       2^{\frac{j q_{1} }{2}} 2^{j s_{1} q_{1} }   |  \beta_{p,j,k} |^{ q_{1} }  |\chi_{j,k}(x)|  \Big )^{\frac{r_{1} }{q_{1}}}  dx \Big )^{\frac{ \theta}{r_{1}}}    .
\end{align*}
Now when $ \varepsilon  $ tends to zero we find
\begin{align*}
\Vert \lambda \vert f^{s}_{r,q}(m)   \Vert & \leq \Vert \lambda \vert  ( f^{s_{0}}_{r_{0},q_{0}}(m) )^{1 - \theta}    ( f^{s_{1}}_{r_{1},q_{1}}(m) )^{ \theta}   \Vert  \Vert \alpha \vert  f^{s_{0}}_{r_{0},q_{0}}(m)   \Vert^{1 - \theta}   \Vert \beta \vert  f^{s_{1}}_{r_{1},q_{1}}(m)   \Vert^{\theta} \\
& \leq  \Vert \lambda \vert  ( f^{s_{0}}_{r_{0},q_{0}}(m) )^{1 - \theta}    ( f^{s_{1}}_{r_{1},q_{1}}(m) )^{ \theta}   \Vert  .
\end{align*}
Here in the last step we used the special choice of the sequences $ \alpha  $ and $ \beta  $. So step 1 of the proof is complete. 

\textit{Step 2.} Now we prove the converse inequality. For that purpose let $  \lambda = \{ c_{p,j,k}   \}_{p,j,k}   $ be a sequence with $  \lambda \in f^{s}_{r,q}(m)   $. Again we proceed like it is described in the proof of Theorem 8.2 from \cite{FJ2}. That means in a first step we assume $ \frac{r_{0}}{q_{0}} \leq  \frac{r_{1}}{q_{1}}  $. Then the reverse case can be done with similar arguments by interchanging $  f^{s_{0}}_{r_{0},q_{0}}(m)  $ and $   f^{s_{1}}_{r_{1},q_{1}}(m)  $. Now for each $ l \in \mathbb{Z}  $ we define sets $ A_{l} \subset \mathbb{R}   $ by
\begin{align*}
A_{l} := \Big \{  x \in \mathbb{R} \ : \   \Big[ \sum_{p = 0}^{\infty}  \sum_{j= -1}^\infty \sum_{k \in \mathbb{Z}}  ( p + 1 )^{  (2m+3)^{2} q}  2^{jsq} 2^{\frac{jq}{2}}    |c_{p,j,k}|^{q}  |\chi_{j,k}(x)|^{q}   \Big]^{\frac{1}{q}}   > 2^{l}    \Big \} .
\end{align*}
Moreover each triple $(p,j,k)$ can be associated with a dyadic interval $ Q_{j,k}  $. Then $ K_{l}    $ is the collection of all triples $(p,j,k)$ that fulfill
\begin{align*}
| Q_{j,k} \cap A_{l} | \geq \frac{|Q_{j,k}|}{2} \qquad \qquad \mbox{and} \qquad \qquad | Q_{j,k} \cap A_{l+1} | < \frac{|Q_{j,k}|}{2}.
\end{align*}
Notice that if for one triple $(p,j,k)$ we observe $ (p,j,k) \not \in  \cup_{l \in \mathbb{Z}} K_{l}      $, then we have $  c_{p,j,k} = 0   $. In what follows we need some additional notation. Let us write
\begin{align*}
\gamma = 1 - \frac{r q_{0}}{q r_{0}} \qquad \qquad \mbox{and} \qquad \qquad \delta =  1 - \frac{r q_{1}}{q r_{1}} .
\end{align*}
Moreover we want to work with numbers $u$ and $v$ given by 
\begin{align*}
u = s + \frac{1}{2} - \frac{q_{0}}{q} \Big (  s_{0} + \frac{1}{2} \Big ) \qquad \qquad \mbox{and} \qquad \qquad v = s + \frac{1}{2} - \frac{q_{1}}{q} \Big (  s_{1} + \frac{1}{2} \Big ).
\end{align*}
We can use the notation we just explained to define for each triple $(p,j,k)\in K_{l} $ numbers $ A_{p,j,k}  $ and $ B_{p,j,k}  $ by
\begin{align*}
 A_{p,j,k} = 2^{l \gamma } |Q_{j,k}|^{u}    \qquad \qquad \mbox{and} \qquad \qquad    B_{p,j,k} =  2^{l \delta } |Q_{j,k}|^{v} .
\end{align*}
Now we define sequences $ \alpha = \{ \alpha_{p,j,k}   \}_{p,j,k} $ and $ \beta = \{ \beta_{p,j,k}   \}_{p,j,k} $ via 
\begin{align*}
\alpha_{p,j,k}  = \Big ( (p+1)^{(1 - \frac{q_{0}}{q})(2m+3)^{2}} \frac{|c_{p,j,k} |}{A_{p,j,k}}  \Big )^{\frac{q}{q_{0}}}  \qquad \mbox{and} \qquad  \beta_{p,j,k}  = \Big ( (p+1)^{(1 - \frac{q_{1}}{q})(2m+3)^{2}} \frac{|c_{p,j,k} |}{B_{p,j,k}}  \Big )^{\frac{q}{q_{1}}}.
\end{align*}
In case that $ (p,j,k) \not \in  \cup_{l \in \mathbb{Z}} K_{l}      $, then we put $  \alpha_{p,j,k} = \beta_{p,j,k} = 0      $. Now we prove that we have 
\begin{equation}\label{prof_caldp_coe}
|  c_{p,j,k} | = |  \alpha_{p,j,k} |^{1 - \theta} | \beta_{p,j,k}  |^{\theta}.
\end{equation}
For the case $(p,j,k) \not \in  \cup_{l \in \mathbb{Z}} K_{l}$ this is obvious because then we have $  c_{p,j,k} = \alpha_{p,j,k} = \beta_{p,j,k} = 0     $. Now let $(p,j,k)\in K_{l} $. Then we have
\begin{align*}
& |  \alpha_{p,j,k} |^{1 - \theta} | \beta_{p,j,k}  |^{\theta} \\
& \qquad = |c_{p,j,k} |^{\frac{q(1 - \theta)}{q_{0}}} 2^{-l \gamma \frac{q(1 - \theta)}{q_{0}} } 2^{ju \frac{q(1 - \theta)}{q_{0}}}  |c_{p,j,k} |^{\frac{q \theta}{q_{1}}} 2^{- l \delta \frac{q \theta}{q_{1}} } 2^{jv \frac{q \theta}{q_{1}}}  (p+1)^{(2m+3)^{2}(\frac{q}{q_{0}}(1 - \theta)+ \frac{q}{q_{1}} \theta  - 1)}  \\
& \qquad = |c_{p,j,k} |^{\frac{q(1 - \theta)}{q_{0}} + \frac{q \theta}{q_{1}}}   2^{-l \gamma \frac{q(1 - \theta)}{q_{0}} - l \delta \frac{q \theta}{q_{1}}  }   2^{ju \frac{q(1 - \theta)}{q_{0}} + jv \frac{q \theta}{q_{1}} }  (p+1)^{(2m+3)^{2}(\frac{q}{q_{0}}(1 - \theta)+ \frac{q}{q_{1}} \theta  - 1)}  . 
\end{align*}
We observe
\begin{align*}
\frac{q(1 - \theta)}{q_{0}} + \frac{q \theta}{q_{1}} = 1 \qquad \mbox{and} \qquad \gamma \frac{q(1 - \theta)}{q_{0}} + \delta \frac{q \theta}{q_{1}} = 0  \qquad \mbox{and} \qquad u \frac{q(1 - \theta)}{q_{0}} + v \frac{q \theta}{q_{1}} = 0 .
\end{align*}
This shows formula \eqref{prof_caldp_coe}. Next we want to prove that we have
\begin{equation}\label{prof_caldp_coe2}
\Vert \alpha \vert f^{s_{0}}_{r_{0},q_{0}}(m)   \Vert \leq C_{1} \Vert  \lambda \vert f^{s}_{r,q}(m)  \Vert^{\frac{r}{r_{0}}}  \qquad \mbox{and}  \qquad \Vert \beta \vert f^{s_{1}}_{r_{1},q_{1}}(m)   \Vert \leq C_{1} \Vert  \lambda \vert f^{s}_{r,q}(m)  \Vert^{\frac{r}{r_{1}}} .
\end{equation}
To see this we follow the strategy explained on page 94 in \cite{FJ2}.  At first we use the definitions to find
\begin{align*}
& \Vert \alpha \vert f^{s_{0}}_{r_{0},q_{0}}(m)   \Vert^{r_{0}} \\
& \qquad \leq C_{2}   \int_{\mathbb{R}} \Big[ \sum_{l \in \mathbb{Z}}   \sum_{p,j,k \in K_{l}}   ( p + 1 )^{  (2m+3)^{2} q}  2^{js_{0}q_{0}} 2^{\frac{jq_{0}}{2} }    A_{p,j,k}^{-q}  |c_{p,j,k}|^{q}  \chi_{j,k}(x)  \chi_{A_{l}}(x) \Big]^{\frac{r_{0}}{q_{0}}}  dx  \\
& \qquad \leq C_{3}   \int_{\mathbb{R}} \Big[ \sum_{l \in \mathbb{Z}}  \chi_{A_{l}}(x)  \sum_{p,j,k \in K_{l}}  2^{- l \gamma q}   ( p + 1 )^{  (2m+3)^{2} q}  2^{jsq} 2^{\frac{jq}{2}}    |c_{p,j,k}|^{q}  \chi_{j,k}(x)   \Big]^{\frac{r_{0}}{q_{0}}}  dx  . 
\end{align*}
Notice that we have $ \gamma \leq 0   $. Moreover on each set $ A_{l} $ we can estimate the number $ 2^{l} $ from above using the definition of $A_{l}$. So we get
\begin{align*}
& \Vert \alpha \vert f^{s_{0}}_{r_{0},q_{0}}(m)   \Vert^{r_{0}} \\
& \qquad \leq C_{4}   \int_{\mathbb{R}} \Big[ \sum_{l \in \mathbb{Z}}  \chi_{A_{l}}(x)   \Big[ \sum_{p = 0}^{\infty}  \sum_{j= -1}^\infty \sum_{k \in \mathbb{Z}}  ( p + 1 )^{  (2m+3)^{2} q}  2^{jsq} 2^{\frac{jq}{2}}    |c_{p,j,k}|^{q}  |\chi_{j,k}(x)|^{q}   \Big]^{- 1 + \frac{r}{q} \frac{q_{0}}{r_{0}}}  \\
& \qquad \qquad \qquad \qquad \qquad  \Big [ \sum_{p',j',k' \in K_{l}}    ( p' + 1 )^{  (2m+3)^{2} q}  2^{j'sq} 2^{\frac{j'q}{2}}    |c_{p',j',k'}|^{q}  \chi_{j',k'}(x) \Big ]  \Big]^{\frac{r_{0}}{q_{0}}}  dx  \\ 
& \qquad \leq C_{5} \Vert  \lambda \vert f^{s}_{r,q}(m)  \Vert^{r} . 
\end{align*}  
This proves the first part of \eqref{prof_caldp_coe2}. For the proof of the second part we can use similar arguments. The requiered modifications are explained on page 95 in \cite{FJ2}. So we omit the details. Next we observe that \eqref{prof_caldp_coe} also can be written as 
\begin{align*}
|  c_{p,j,k} |   =  \Vert  \lambda \vert f^{s}_{r,q}(m)  \Vert   \Big ( \frac{ | \alpha_{p,j,k} |}{\Vert  \lambda \vert f^{s}_{r,q}(m)  \Vert^{\frac{r}{r_{0}}}} \Big )^{1 - \theta}  \Big ( \frac{ | \beta_{p,j,k}  |}{\Vert  \lambda \vert f^{s}_{r,q}(m)  \Vert^{\frac{r}{r_{1}}}} \Big )^{\theta} .
\end{align*} 
Now we can use Definition \ref{def_cald_P} and the estimates from \eqref{prof_caldp_coe2} to find
\begin{align*}
 \Vert \lambda \vert  ( f^{s_{0}}_{r_{0},q_{0}}(m) )^{1 - \theta}    ( f^{s_{1}}_{r_{1},q_{1}}(m) )^{ \theta}   \Vert &   \leq \Vert \alpha \vert f^{s_{0}}_{r_{0},q_{0}}(m)   \Vert^{1 - \theta}  \Vert \beta \vert f^{s_{1}}_{r_{1},q_{1}}(m)   \Vert^{\theta} \\
&   \leq C_{6}  \Vert  \lambda \vert f^{s}_{r,q}(m)  \Vert^{\frac{r}{r_{0}}(1 - \theta)} \Vert  \lambda \vert f^{s}_{r,q}(m)  \Vert^{\frac{r}{r_{1}} \theta}  \\
&   = C_{6}  \Vert  \lambda \vert f^{s}_{r,q}(m)  \Vert  .
\end{align*}
This in combination with \eqref{prof_caldp_coe} and $   \lambda \in f^{s}_{r,q}(m)  $ yields the desired result. The proof is complete.  
\end{proof}

It turns out that especially in the case of sequence spaces Calderon products $  X_{0}^{1 - \theta}  X_{1}^{ \theta}   $ are very strongly connected with the interpolation spaces $  [X_0,X_1]_\theta  $ that one obtains by applying Calderons first complex interpolation method. Therefore with a view to Lemma \ref{lem_cald_seq} it is not surprising that we have the following interpolation result. 

\begin{lem}\label{lem_calfirstcom_seq}
Let $ s_{0}, s_{1} > 0  $, $ 0 < r_{0}, r_{1} < \infty $ and $  0 < q_{0},q_{1} < \infty $. Let $ m \in \mathbb{N}  $. Let $ 0 < \theta < 1 $ and put 
\begin{equation}\label{comp_standcon}
s = ( 1 - \theta ) s_{0} + \theta s_{1} , \qquad  \qquad \frac{1}{r} = \frac{1 - \theta}{r_{0}} + \frac{\theta}{r_{1}} , \qquad  \qquad \frac{1}{q} = \frac{1 - \theta}{q_{0}} + \frac{\theta}{q_{1}} .
\end{equation}
Then for Calderons first complex interpolation method we observe
\begin{align*}
 [ f^{s_{0}}_{r_{0},q_{0}}(m) ,  f^{s_{1}}_{r_{1},q_{1}}(m) ]_{\theta} =  f^{s}_{r,q}(m)  .
\end{align*}
\end{lem} 

\begin{proof}
Lemma \ref{lem_calfirstcom_seq} is an easy consequence of Lemma \ref{lem_cald_seq} and the formula 
\begin{equation}\label{calp=calcom1}
[ f^{s_{0}}_{r_{0},q_{0}}(m) ,  f^{s_{1}}_{r_{1},q_{1}}(m) ]_{\theta} = ( f^{s_{0}}_{r_{0},q_{0}}(m) )^{1 - \theta}    ( f^{s_{1}}_{r_{1},q_{1}}(m) )^{ \theta}.
\end{equation}
More general versions of formula \eqref{calp=calcom1} already can be found in the literature. For the case of Banach spaces it can be found on page 125 in \cite{ca64}. We also refer to the remark before Corollary 8.3 in \cite{FJ2}. For quasi-Banach spaces \eqref{calp=calcom1} follows from Theorem 3.4 in \cite{KM1998}, see also Theorem 7.9 in \cite{KMM}. We only have to check that all conditions that can be found in Theorem 3.4 in \cite{KM1998} are fulfilled in our case. So on the one hand like it is described in the remark after Theorem 3.4, we have to verify the following property. Let $  \lambda = \{ c_{p,j,k}   \}_{p,j,k} \in f^{s}_{r,q}(m)   $ be a sequence and let $\{ \lambda^{n} \}_{n} = \{ \{ c_{p,j,k}^{n}   \}_{p,j,k} \}_{n}   $ be a sequence of sequences such that for all $ n \in \mathbb{N}   $ we have $ | c_{p,j,k}^{n} | \leq | c_{p,j,k} |     $ and $ \lim_{n \rightarrow \infty } c_{p,j,k}^{n} = c_{p,j,k}  $ for all $ p,j,k  $. Then we need that we also have
\begin{equation}\label{calp=calcom1_eq1}
\lim_{n \rightarrow \infty } \Vert \lambda^{n} - \lambda   \vert   f^{s}_{r,q}(m)  \Vert = 0.
\end{equation}
But this follows from the dominated convergence theorem (like it was also observed in the proof of Corollary 8.3 in \cite{FJ2}). Second we have to prove that the spaces $ f^{s}_{r,q}(m)  $ are analytically convex. The term    analytically convex for example is explained on page 21 in \cite{KMM}, see formula (7.7). But the fact that the sequence spaces $ f^{s}_{r,q}(m)  $ are analytically convex can be proved in the same way as for the original sequence spaces without polynomial enrichment. For that purpose we refer to Lemma 7.6 and Proposition 7.7 in \cite{KMM}. So all conditions that can be found in Theorem 3.4 in \cite{KM1998} are fulfilled. The proof is complete.
\end{proof}


Now we are prepared to prove a counterpart of Proposition \ref{res_Fspq_m-1} that also holds for $ s \geq m - 1  $. Here we obtain the following result. 

\begin{prop}\label{res_Fspq_finalup}
Let $  0 < r < \infty $, $ 0 < q < \infty  $ and $ m \in \mathbb{N}  $ with $ m \geq 2  $. Let $ s \geq m - 1   $. Assume that either condition $ (II)$ or condition $(III)$ from Theorem \ref{mainresult1} is fulfilled. Let $ f \in \mathcal{S}'(\mathbb{R})  $ such that there exists a representation 
\begin{equation}\label{eq-inf4}
f = \sum_{p \geq 0} \sum_{j \geq -1} \sum_{k \in \mathbb{Z}} c_{p,j,k} \psi_{p,j,k}
\end{equation}
with convergence in $ \mathcal{S}'(\mathbb{R})  $. Then there exists a $ C > 0  $ independent from $ f $ such that  
\begin{align*}
& \|\, f \, |F^s_{r,q}(\mathbb{R})\| \leq C \inf_{\eqref{eq-inf4}} \Big\| \Big[ \sum_{p = 0}^{\infty}  \sum_{j= -1}^\infty \sum_{k \in \mathbb{Z}}  ( p + 1 )^{(2m+3)^{2} q}  2^{jsq} 2^{\frac{jq}{2}}    |c_{p,j,k}|^{q}  |\chi_{j,k}(x)|^{q}   \Big]^{\frac{1}{q}}\Big|  L_{r} (\mathbb{R}) \Big\| .
\end{align*}
Here the infimum is taken over all sequences $ \{ c_{p,j,k}  \}_{p,j,k} \subset \mathbb{C} $ such that \eqref{eq-inf4} is fulfilled.
\end{prop}

\begin{proof}
To prove this result we use a combination of Proposition \ref{stab_r=q_full}, Proposition \ref{res_Fspq_m-1} and complex interpolation. From Proposition \ref{stab_r=q_full} we know
\begin{align*}
 \|\, f \, |F^{s_{1}}_{r_{1},r_{1}}(\mathbb{R})\| \leq C \inf_{\eqref{eq-inf4}} \Vert \lambda \vert f^{s_{1}}_{r_{1},r_{1}}(m)   \Vert  
\end{align*} 
for
\begin{equation}\label{inter-eq2}
0 < r_{1} = q_{1} < \infty  , \quad  \max \Big(0, \frac {1}{r_{1}} -1  \Big) < s_{1} <  \min \Big ( m - 1 + \frac{1}{r_{1}} , m  \Big ) .
\end{equation}  
On the other hand from Proposition \ref{res_Fspq_m-1} we know that we have
\begin{align*}
 \|\, f \, |F^{s_{0}}_{r_{0},q_{0}}(\mathbb{R})\| \leq C \inf_{\eqref{eq-inf4}} \Vert \lambda \vert f^{s_{0}}_{r_{0},q_{0}}(m)   \Vert  
\end{align*} 
for
\begin{equation}\label{inter-eq1}
0 < r_{0} < \infty , \quad  0 < q_{0} < \infty , \quad \max \Big(0, \frac {1}{r_{0}} -1 , \frac {1}{q_{0}} -1  \Big) < s_{0}  <  m - 1 .
\end{equation}
Now we want to use complex interpolation, see Lemma \ref{lem_comint_Fsrq1} and Lemma \ref{lem_calfirstcom_seq}. We use some ideas from the proof of Proposition 2.6 in \cite{Tr10}. We show that we can choose $ 0 < \theta < 1  $ and $ (s_{0},r_{0},q_{0})    $ as well as $  (s_{1},r_{1},q_{1})   $ in such a way that each triple $  (s,r,q)   $ that fulfills either condition $ (II)$ or condition $(III)$ from Theorem \ref{mainresult1} can be reached via complex interpolation, see \eqref{comp_standcon}. For that purpose let $  (s,r,q)   $ be fixed, such that $ (II)$ or $(III)$ are fulfilled. Now we put $  r = r_{0} = r_{1}  = q_{1}    $. Furthermore we define $ s_{0} = m - 1 - \varepsilon   $ for some small $ \varepsilon > 0  $ and 
\begin{align*}
s_{1} = \min \Big ( m - 1 + \frac{1}{r}   , m  \Big ) - \varepsilon .
\end{align*}
Then of course we have $ \frac{1}{r} = \frac{1 - \theta}{r_{0}} + \frac{\theta}{r_{1}} $ for all $ 0 < \theta < 1    $. Moreover now $ \theta = \theta(s,r,q,m)    $ is defined via $ s = ( 1 - \theta ) s_{0} + \theta s_{1} $ which results in
\begin{align*}
\theta = \frac{s - m + 1 + \varepsilon}{ 1 + \min(  \frac{1}{r} - 1 ,    0   )} .
\end{align*}
Now in the setting of complex interpolation we have the relation $ \frac{1}{q} = \frac{1 - \theta}{q_{0}} + \frac{\theta}{r} $. Recall that the desired result holds for all $ q_{0}  $ such that $ \frac{1}{m} < q_{0} < \infty   $. With $ q_{0} \rightarrow \infty   $ and $ q_{0} \rightarrow \frac{1}{m}    $ we find
\begin{align*}
\frac{1}{q} \in \Big ( \frac{s-m+1 }{r (1 + \min (\frac{1}{r} - 1 , 0))} , \frac{r(1 + \min (\frac{1}{r} - 1, 0)) - r(s-m+1) + \frac{1}{m} (s-m+1)}{\frac{r}{m} (1 + \min (\frac{1}{r} - 1, 0)) }  \Big ) .
\end{align*}
In the case $ r \geq 1  $ this implies
\begin{align*}
s < m - 1 + \frac{1}{q} \qquad \qquad \mbox{and} \qquad \qquad \frac{1}{q} < s + 1 - mr(s + 1 - m) .
\end{align*}
Here the first condition seems to be natural. The second condition sanctions large $s$. For $ s = m -1  $ it reads as $ \frac{1}{q} < m   $ which is the best we can reach. In the worst case for $ s  $ near to $ \min  ( m - 1 + \frac{1}{r} , m - 1 + \frac{1}{q}    )  $ we get $ \frac{1}{q} < \frac{1}{r}  $. For $ r < 1   $ we get
\begin{align*}
s < m - 1 + \frac{r}{q} \qquad \qquad \mbox{and} \qquad \qquad \frac{1}{q} < - sm + m^{2} + \frac{s}{r} - \frac{m}{r} + \frac{1}{r} .
\end{align*}
Here the first condition only plays a role for $ r < q  $. The second condition again sanctions large $ s $. For $ s = m -1  $ it becomes $ \frac{1}{q} < m   $. For $ s $ near $ m $ we obtain $ \frac{1}{q} < \frac{1}{r}  $. Now we consider the linear operator
\begin{align*}
J : \qquad \qquad  \lambda = \{ c_{p,j,k}   \}_{p,j,k}  \qquad \mapsto \qquad f = \sum_{p \geq 0} \sum_{j \geq -1} \sum_{k \in \mathbb{Z}} c_{p,j,k} \psi_{p,j,k} .
\end{align*}
We recall that the complex interpolation method $ [ \cdot , \cdot   ]_{\theta}   $ has the so-called interpolation property for linear operators. For the Banach space case this is known since many years. For quasi-Banach spaces we refer to \cite{KMM}, see Theorem 8.1. One may also consult the remark before Theorem 1.22. in \cite{Tr10}. Now we use the definition of the operator norm in combination with the interpolation property and the Lemmas \ref{lem_comint_Fsrq1} and \ref{lem_calfirstcom_seq} to find
\begin{align*}
\sup_{\lambda \not = 0} \frac{\Vert f \vert F^{s}_{r,q}(\mathbb{R}) \Vert}{\Vert \lambda \vert f^{s}_{r,q}(m)   \Vert  } &  = \Vert J \vert f^{s}_{r,q}(m) \rightarrow F^{s}_{r,q}(\mathbb{R})   \Vert \\
& = \Vert J \vert   [ f^{s_{0}}_{r_{0},q_{0}}(m) ,  f^{s_{1}}_{r_{1},q_{1}}(m)  ]_{\theta}    \rightarrow [ F^{s_{0}}_{r_{0},q_{0}}(\mathbb{R}) ,  F^{s_{1}}_{r_{1},q_{1}}(\mathbb{R})   ]_{\theta}   \Vert    \\
& \leq \Vert J \vert    f^{s_{0}}_{r_{0},q_{0}}(m)    \rightarrow  F^{s_{0}}_{r_{0},q_{0}}(\mathbb{R})    \Vert^{1 - \theta}  \Vert J \vert    f^{s_{1}}_{r_{1},q_{1}}(m)    \rightarrow  F^{s_{1}}_{r_{1},q_{1}}(\mathbb{R})    \Vert^{\theta} \\
& =  \sup_{\lambda \not = 0} \frac{\Vert f \vert F^{s_{0}}_{r_{0},q_{0}}(\mathbb{R}) \Vert^{1 - \theta}}{\Vert \lambda \vert f^{s_{0}}_{r_{0},q_{0}}(m)   \Vert^{1 - \theta}  }    \sup_{\lambda \not = 0} \frac{\Vert f \vert F^{s_{1}}_{r_{1},q_{1}}(\mathbb{R}) \Vert^{ \theta}}{\Vert \lambda \vert f^{s_{1}}_{r_{1},q_{1}}(m)   \Vert^{ \theta}  }  .
\end{align*}
Next we apply the results from Proposition \ref{res_Fspq_m-1} and Proposition \ref{stab_r=q_full}. Then we find
\begin{align*}
\sup_{\lambda \not = 0} \frac{\Vert f \vert F^{s}_{r,q}(\mathbb{R}) \Vert}{\Vert \lambda \vert f^{s}_{r,q}(m)   \Vert  } \leq C  \sup_{\lambda \not = 0} \frac{\Vert \lambda \vert f^{s_{0}}_{r_{0},q_{0}}(m)   \Vert^{1 - \theta}}{\Vert \lambda \vert f^{s_{0}}_{r_{0},q_{0}}(m)   \Vert^{1 - \theta}  }    \sup_{\lambda \not = 0} \frac{\Vert \lambda \vert f^{s_{1}}_{r_{1},q_{1}}(m)   \Vert^{ \theta} }{\Vert \lambda \vert f^{s_{1}}_{r_{1},q_{1}}(m)   \Vert^{ \theta}  } = C  .
\end{align*}
Therefore for all $ \lambda \in  f^{s}_{r,q}(m)   $ we get $ \Vert f \vert F^{s}_{r,q}(\mathbb{R}) \Vert  \leq C \Vert \lambda \vert f^{s}_{r,q}(m)   \Vert $. The proof is complete. 
\end{proof}

\subsection{Lower estimates}

In this section we want to supplement the Propositions \ref{res_Fspq_m-1} and \ref{res_Fspq_finalup} by proving the corresponding lower estimates. For that purpose in a first step we work with $ p = 0$, namely the case that the quarklets do not have any polynomial enrichment. Recall that $ \tilde{\psi}_{j,k} $ are defined as in \eqref{def_biort_wav1} and \eqref{def_biort_wav2}.

\begin{prop}\label{thm_spline_trCW}
Let $ m \in \mathbb{N}  $. Let $ 0 < r < \infty   $ and $ 0 < q \leq \infty   $ and $ \sigma_{r,q}  < s < m $. Let $ f \in F^{s}_{r,q}(\mathbb{R})  $. For $ j \geq -1  $ and $ k \in \mathbb{Z}   $ we write
\begin{equation}\label{lambdajk}
\lambda_{j,k} =  \int_{- \infty}^{\infty} f(x) \tilde{\psi}_{j,k}(x) dx .
\end{equation}
Then there is a constant $ C > 0  $ independent of $ f \in F^{s}_{r,q}(\mathbb{R})  $ such that
\begin{equation}\label{thm_spline_trCW_eq1}
\Big \Vert \Big ( \sum_{j = -1}^{\infty} \sum_{k \in \mathbb{Z}} 2^{jsq} 2^{\frac{jq}{2}} |\lambda_{j,k}|^{q} |\chi_{j,k}(x)|^{q}    \Big )^{\frac{1}{q}} \Big \vert  L_{r}(\mathbb{R})  \Big \Vert \leq C  \|\, f \, |F^s_{r,q}(\mathbb{R})\| .
\end{equation}
\end{prop}

\begin{proof}
In general all what we need for the proof is known and can be found in the proof of Theorem 2.49 in \cite{Tr10}. Therefore we will be rather brief in what follows. The main idea for the proof is to interpret the functions $ 2^{\frac{j}{2}} \tilde{\psi}_{j,k} $ as kernels of local means according to Definition 1.9 in \cite{Tr10}. We use Definition 1.9 with $ A = 0   $ and $ B = m  $. The numbers $ 2^{\frac{j}{2}} \lambda_{j,k}  $ can be seen as local means according to Definition 1.13 in \cite{Tr10}. Therefore we can apply Theorem 1.15 in \cite{Tr10}, see also Theorem 15 in \cite{Tr08}. So for $ \sigma_{r,q}  < s < m $ and for all $ f \in F^{s}_{r,q}(\mathbb{R})   $ there is a constant $ C $ independent of $ f \in F^{s}_{r,q}(\mathbb{R})   $ such that we have \eqref{thm_spline_trCW_eq1}. 
\end{proof}

Under some additional conditions on the parameters it is possible to improve Proposition \ref{thm_spline_trCW}. So we observe the following.

\begin{prop}\label{Fspq_prefinallow}
Let $ s \in \mathbb{R}  $, $  0 < r < \infty $, $ 0 < q < \infty  $ and $ m \in \mathbb{N}  $ with $ m \geq 2  $. Assume that the parameters satisfy one of the conditions $ (I)  $, $(II)  $ or $  (III)  $ from Theorem \ref{mainresult1}. Let $ f \in F^{s}_{r,q}(\mathbb{R})  $. For $ j \geq -1  $ and $ k \in \mathbb{Z}   $ the numbers $\lambda_{j,k} $ are defined as in \eqref{lambdajk}. Then $ f $ can be represented as
\begin{equation}\label{rep_ohnequa}
f = \sum_{j = -1}^{\infty} \sum_{k \in \mathbb{Z}} \lambda_{j,k} \psi_{j,k}
\end{equation}
with convergence in $ \mathcal{S}'(\mathbb{R})  $. Moreover there is a constant $ C > 0  $ independent of $ f \in F^{s}_{r,q}(\mathbb{R})  $ such that \eqref{thm_spline_trCW_eq1} is fulfilled.
\end{prop}

\begin{proof}
Most of the result already has been proved after Proposition \ref{thm_spline_trCW}. Therefore it remains to show the representation in \eqref{rep_ohnequa}. For that purpose let $ f \in F^{s}_{r,q}(\mathbb{R})  $. We are interested in the function
\begin{align*}
g = \sum_{j = -1}^{\infty} \sum_{k \in \mathbb{Z}} \lambda_{j,k} \psi_{j,k} .
\end{align*}
It has the same shape as the functions we considered in the Propositions \ref{res_Fspq_m-1} and \ref{res_Fspq_finalup} when we put $ c_{p,j,k} = 0   $ for $ p > 0  $ in \eqref{eq-inf3}. Moreover all conditions on the parameters that are given there are fulfilled. Therefore we can use the calculations from there to find
\begin{align*}
\|\, g \, |F^s_{r,q}(\mathbb{R})\| \leq C_{1}   \Big \Vert \Big ( \sum_{j = -1}^{\infty} \sum_{k \in \mathbb{Z}} 2^{jsq} 2^{\frac{jq}{2}} |\lambda_{j,k}|^{q} |\chi_{j,k}(x)|^{q}    \Big )^{\frac{1}{q}} \Big \vert  L_{r}(\mathbb{R})  \Big \Vert \leq C_{2}  \|\, f \, |F^s_{r,q}(\mathbb{R})\| < \infty ,
\end{align*}
where we also used Proposition \ref{thm_spline_trCW} and $ f \in F^{s}_{r,q}(\mathbb{R})   $.  With other words we have $ g \in F^s_{r,q}(\mathbb{R})   $. Now we want to prove that we have $ f = g  $ in the sense of $ \mathcal{S}'(\mathbb{R})  $. For that purpose we show that for all $ \eta \in \mathcal{S}(\mathbb{R})  $ we have $ \langle f - g ,  \eta \rangle_{L_{2}(\mathbb{R})} = 0 $. Since $  \mathcal{S}(\mathbb{R})  \subset L_{2}(\mathbb{R})   $ we can write
\begin{align*}
\eta = \sum_{j = -1}^{\infty}  \sum_{k \in \mathbb{Z}} \left\langle \eta , \psi_{j,k}     \right\rangle_{L_{2}(\mathbb{R})} \tilde{\psi}_{j,k}
\end{align*}
with convergence in $ L_{2}(\mathbb{R})  $, see \eqref{def_biort_wav3}. From convergence in $ L_{2}(\mathbb{R})  $ it follows the convergence pointwise almost everywhere for an appropriate subsequence. Now let $ j' \geq -1 $ and $ k' \in \mathbb{Z}   $ be fixed. Then we find
\begin{align*}
& \int_{\mathbb{R}} (f(x) - g(x)) \langle \eta , \psi_{j',k'}     \rangle_{L_{2}(\mathbb{R})} \tilde{\psi}_{j',k'}(x) dx \\
& \qquad =  \int_{\mathbb{R}} \Big (f(x) - \sum_{j = -1}^{\infty} \sum_{k \in \mathbb{Z}} \lambda_{j,k} \psi_{j,k}(x) \Big) \langle \eta , \psi_{j',k'}     \rangle_{L_{2}(\mathbb{R})} \tilde{\psi}_{j',k'}(x) dx \\
& \qquad = \langle \eta , \psi_{j',k'}     \rangle_{L_{2}(\mathbb{R})}    \langle f , \tilde{\psi}_{j',k'} \rangle_{L_{2}(\mathbb{R})}   - \langle \eta , \psi_{j',k'}     \rangle_{L_{2}(\mathbb{R})}   \sum_{j = -1}^{\infty} \sum_{k \in \mathbb{Z}} \lambda_{j,k}  \int_{\mathbb{R}} \psi_{j,k}(x) \tilde{\psi}_{j',k'}(x)    dx .
\end{align*}
Now we use the definition of the numbers $ \lambda_{j,k}  $ and the fact that the wavelet bases associated to $ \psi  $ and $ \tilde{\psi}  $ are biorthogonal, see \eqref{biorto1} and \eqref{biorto2}. Therefore we obtain
\begin{align*}
& \int_{\mathbb{R}} (f(x) - g(x)) \langle \eta , \psi_{j',k'}     \rangle_{L_{2}(\mathbb{R})} \tilde{\psi}_{j',k'}(x) dx \\
& \qquad = \langle \eta , \psi_{j',k'}     \rangle_{L_{2}(\mathbb{R})}    \langle f , \tilde{\psi}_{j',k'} \rangle_{L_{2}(\mathbb{R})}   - \langle \eta , \psi_{j',k'}     \rangle_{L_{2}(\mathbb{R})} \langle f , \tilde{\psi}_{j',k'} \rangle_{L_{2}(\mathbb{R})} = 0 .  
\end{align*}
This calculation can be extended to every linear combination of $  \left\langle \eta , \psi_{j,k}     \right\rangle_{L_{2}(\mathbb{R})} \tilde{\psi}_{j,k} $. This observation in combination with our considerations concerning convergence yields $ f = g   $ in the sense of $ \mathcal{S}'(\mathbb{R})  $. So the proof is complete. 
\end{proof}

Now we are well prepared to prove the following lower estimate for the full quarklet system. 

\begin{prop}\label{Fspq_finallow}
Let $ s \in \mathbb{R}  $, $  0 < r < \infty $, $ 0 < q < \infty  $ and $ m \in \mathbb{N}  $ with $ m \geq 2  $. Moreover the parameters satisfy one of the conditions $ (I)  $, $(II)  $ or $  (III)  $ from Theorem \ref{mainresult1}. Let $ f \in F^{s}_{r,q}(\mathbb{R})  $.  Then there exists a sequence $ \{ c_{p,j,k}  \}_{p \in \mathbb{N}_{0}, k \in \mathbb{Z}, j \geq -1}   $ such that $ f $ can be represented as
\begin{equation}\label{rep_qua}
f = \sum_{p \geq 0} \sum_{j = -1}^{\infty} \sum_{k \in \mathbb{Z}} c_{p,j,k} \psi_{p,j,k}
\end{equation}
with convergence in $ \mathcal{S}'(\mathbb{R})  $. Moreover there is a constant $ C > 0  $ independent of $ f \in F^{s}_{r,q}(\mathbb{R})  $ such that
\begin{align*}
\inf_{\eqref{rep_qua}} \Big\| \Big[ \sum_{p = 0}^{\infty}  \sum_{j= -1}^\infty \sum_{k \in \mathbb{Z}}  ( p + 1 )^{(2m+3)^{2} q} 2^{jsq} 2^{\frac{jq}{2}}    |c_{p,j,k}|^{q}  |\chi_{j,k}(x)|^{q}   \Big]^{\frac{1}{q}}\Big|  L_{r} (\mathbb{R}) \Big\|  \leq C  \|\, f \, |F^s_{r,q}(\mathbb{R})\| .
\end{align*}
Here the infimum is taken over all sequences $ \{ c_{p,j,k}  \}   $ such that \eqref{rep_qua} is fulfilled.
\end{prop}

\begin{proof}
Let $ f \in F^{s}_{r,q}(\mathbb{R})  $. Then from Proposition \ref{Fspq_prefinallow} we know that we can write
\begin{equation}\label{pr_swoq}
f = \sum_{j = -1}^{\infty} \sum_{k \in \mathbb{Z}} \lambda_{j,k} \psi_{0,j,k}
\end{equation}
where the numbers $ \lambda_{j,k}  $ are as in \eqref{lambdajk}. The convergence is in $  \mathcal{S}'(\mathbb{R})  $. Now for $ j \in \mathbb{N}_{0} \cup \{ -1 \}  $ and $ k \in \mathbb{Z}   $ we put $ c_{0,j,k} =  \lambda_{j,k}  $. For $ p \in \mathbb{N}  $ we write $ c_{p,j,k} =  0  $. Using this notation we also find
\begin{equation}\label{pr_swq}
f = \sum_{p \in \mathbb{N}_{0}} \sum_{j = -1}^{\infty} \sum_{k \in \mathbb{Z}} c_{p,j,k} \psi_{p,j,k}
\end{equation}
with convergence in $  \mathcal{S}'(\mathbb{R})  $. From Proposition \ref{Fspq_prefinallow} we can also conclude that we have the estimate
\begin{align*}
 \|\, f \, |F^s_{r,q}(\mathbb{R})\|  \geq C   \Big \Vert \Big ( \sum_{j = -1}^{\infty} \sum_{k \in \mathbb{Z}} 2^{jsq} 2^{\frac{jq}{2}} |\lambda_{j,k}|^{q} |\chi_{j,k}(x)|^{q}    \Big )^{\frac{1}{q}} \Big \vert  L_{r}(\mathbb{R})  \Big \Vert  
\end{align*}
with the $  \lambda_{j,k}   $ from \eqref{pr_swoq}. Notice that for $ p = 0  $ we have $ (p + 1)^{(2m+3)^{2} q} =   1^{(2m+3)^{2} q} = 1  $. With this simple observation finally we also find
\begin{align*}
 \|\, f \, |F^s_{r,q}(\mathbb{R})\|  \geq C  \inf_{\eqref{pr_swq}} \Big\| \Big[ \sum_{p = 0}^{\infty}  \sum_{j= -1}^\infty \sum_{k \in \mathbb{Z}}  ( p + 1 )^{(2m+3)^{2} q}  2^{jsq} 2^{\frac{jq}{2}}    |c_{p,j,k}|^{q}  |\chi_{j,k}(x)|^{q}   \Big]^{\frac{1}{q}}\Big|  L_{r} (\mathbb{R}) \Big\| .
\end{align*}
Here the infimum is taken over all sequences $ \{ c_{p,j,k}  \}   $ such that \eqref{pr_swq} is fulfilled. The idea is that now we take the infimum over a larger set of functions which makes it possible to obtain a smaller value for the infimum. So the proof is complete.
\end{proof}

\begin{rem}
There also exist characterizations of the Triebel-Lizorkin spaces $ F^{s}_{r,q}(\mathbb{R})  $ in terms of orthogonal spline wavelets. For that we refer to \cite{Tr10}. Here at the beginning of Chapter 2.5.1. orthogonal spline bases are investigated. In the Theorems 2.46 and 2.49 in \cite{Tr10} they are used to formulate equivalent characterizations for the spaces $ F^{s}_{r,q}(\mathbb{R})  $. Notice that the conditions concerning the parameters that can be found in Theorem 2.49 show many similarities with those we stated in our Theorem \ref{mainresult1}.
\end{rem}

\begin{rem}
It is also possible to describe the Triebel-Lizorkin spaces $ F^{s}_{r,q}(\mathbb{R})   $ in terms of biorthogonal compactly supported Chui-Wang wavelets. For a definition of those wavelets we refer to Theorem 1 in \cite{ChWa}. In \cite{DeUl} the biorthogonal Chui-Wang wavelets are used to prove equivalent quasi-norms for $ F^{s}_{r,q}(\mathbb{R})   $, see Theorem 5.1 and also Theorem 6.2 for higher order Chui-Wang wavelets.  
\end{rem}

\begin{rem}
There exists another version of quarks explained by Triebel, see Section I.2 in \cite{Tr01} and Chapter 1.6 in \cite{Tr06}. These quarks are defined in a different way without using biorthogonal compactly supported B-spline wavelets, see Definition I.2.4 in \cite{Tr01}. Instead they are smoother and more close to atoms. Also the quarks defined in \cite{Tr01} and \cite{Tr06} can be used to describe the Triebel-Lizorkin spaces, see Theorem 1.39 in \cite{Tr06}. Because of the different definition here other conditions concerning the parameters show up. 
\end{rem}

\section{Quarklet characterizations for Triebel-Lizorkin Spaces with negative smoothness}\label{sec_quas<0}

In this section we want to prove quarklet characterizations for Triebel-Lizorkin spaces $ F^{s}_{r,q}(\mathbb{R})   $ with negative smoothness. That means we are interested in the case $ s < 0  $. The strategy to obtain such characterizations is to use the results we already obtained for positive smoothness in combination with some duality arguments. For that purpose we need the dual spaces of $ F^{s}_{r,q}(\mathbb{R})  $ and $  f^{s}_{r,q}(m) $. For $ 1 < r < \infty   $ and $ 1 < q < \infty   $ we define numbers $ r'  $ and $ q'  $  such that
\begin{equation}\label{dual_numb}
\frac{1}{r} + \frac{1}{r'} = 1 \qquad \qquad \mbox{and} \qquad \qquad \frac{1}{q} + \frac{1}{q'} = 1 .
\end{equation}
Using this notation for the Triebel-Lizorkin spaces we know the following. 
\begin{lem}\label{Fsrq_dual}
Let $ s \in \mathbb{R}   $, $ 1 < r < \infty   $ and $ 1 < q < \infty   $. Then we have $(  F^{s}_{r,q}(\mathbb{R})    )' = F^{-s}_{r',q'}(\mathbb{R}) $.
\end{lem}
This result can be found in \cite{Tr83}, see the theorem in Chapter 2.11.2. Now we prove the counterpart for the sequence spaces $ f^{s}_{r,q}(m)  $.

\begin{lem}\label{fsrqm_dual}
Let $ 1 < r < \infty   $, $ 1 < q < \infty   $ and $ s \in \mathbb{R}   $ with $ s \not = 0  $. Let $ m \in \mathbb{N}  $. Then we have 
\begin{align*}
(  f^{s}_{r,q}(m)    )' = f^{-s}_{r',q'}(m) .
\end{align*}
\end{lem}

\begin{proof}
We write down the following proof for $ s>0$. Then the case $ s<0$ can be dealt with the same arguments, but some signs need to be changed. 

\textit{Step 1.} At first we prove $ (  f^{s}_{r,q}(m)    )' \subset f^{-s}_{r',q'}(m)   $. For that purpose we show that each linear functional $ g \in (  f^{s}_{r,q}(m)    )'    $ can be written as
\begin{equation}\label{seq_lin_fun}
g(\lambda) =  \sum_{p = 0}^{\infty}  \sum_{j= -1}^\infty \sum_{k \in \mathbb{Z}} g_{p,j,k} c_{p,j,k}
\end{equation}
for every $ \lambda = \{c_{p,j,k} \} \in  f^{s}_{r,q}(m)  $, where $ \{g_{p,j,k} \}     $ is a sequence with 
\begin{align*}
\Vert g \vert  (  f^{s}_{r,q}(m)    )'   \Vert =   \Big \Vert \Big (  \sum_{p = 0}^{\infty}  \sum_{j= -1}^\infty \sum_{k \in \mathbb{Z}}   ( p + 1 )^{- q' (2m+3)^{2} }  2^{-jsq'} 2^{ \frac{jq'}{2}} | g_{p,j,k} |^{q'}  \chi_{j,k}(x)  \Big )^{\frac{1}{q'}} \Big \vert L_{r'}(\mathbb{R}) \Big \Vert < \infty .
\end{align*}
To see this at first we observe that we can write
\begin{align*}
\Vert \lambda \vert f^{s}_{r,q}(m)   \Vert =   \Big\| \Big[ \sum_{p = 0}^{\infty}  \sum_{j= -1}^\infty \sum_{k \in \mathbb{Z}}     |a_{p,j,k}(x)|^{q}     \Big]^{\frac{1}{q}}\Big|  L_{r} (\mathbb{R}) \Big\| 
\end{align*}
with $  a_{p,j,k}(x) = ( p + 1 )^{(2m+3)^{2} }  2^{js} 2^{\frac{j}{2}}    |c_{p,j,k}|  |\chi_{j,k}(x)| $. Therefore we can use the Proposition in Chapter 2.11.1 in \cite{Tr83}. It tells us that we have a representation
\begin{align*}
g(\lambda) =  \sum_{p = 0}^{\infty}  \sum_{j= -1}^\infty \sum_{k \in \mathbb{Z}} \int_{\mathbb{R}} \bar{g}_{p,j,k}(x)  ( p + 1 )^{(2m+3)^{2} }  2^{js} 2^{\frac{j}{2}}    c_{p,j,k}  |\chi_{j,k}(x)|  dx .
\end{align*}
Here $ \{ \bar{g}_{p,j,k}(x) \}   $  is a sequence of measurable functions with
\begin{align*}
\Vert g \vert  (  f^{s}_{r,q}(m)    )'   \Vert = \Big\| \Big[ \sum_{p = 0}^{\infty}  \sum_{j= -1}^\infty \sum_{k \in \mathbb{Z}}     |\bar{g}_{p,j,k}(x)|^{q'}     \Big]^{\frac{1}{q'}}\Big|  L_{r'} (\mathbb{R}) \Big\| < \infty .
\end{align*}
Next notice that we also can write
\begin{align*}
g(\lambda) & =  \sum_{p = 0}^{\infty}  \sum_{j= -1}^\infty \sum_{k \in \mathbb{Z}}   ( p + 1 )^{(2m+3)^{2} }  2^{js} 2^{\frac{j}{2}}    c_{p,j,k} \int_{Q_{j,k}} \bar{g}_{p,j,k}(x)  dx .
\end{align*}
Consequently we can choose a sequence  $ \{ \tilde{g}_{p,j,k}  \} \subset \mathbb{C}    $ with $ \tilde{g}_{p,j,k} =  \int_{Q_{j,k}} \bar{g}_{p,j,k}(x)  dx   $ such that
\begin{align*}
g(\lambda) & =  \sum_{p = 0}^{\infty}  \sum_{j= -1}^\infty \sum_{k \in \mathbb{Z}}   ( p + 1 )^{(2m+3)^{2} }  2^{js} 2^{\frac{j}{2}}    c_{p,j,k} \tilde{g}_{p,j,k} .
\end{align*}
Therefore it is also possible to choose the functions $ \bar{g}_{p,j,k}  $ to be constant on $ Q_{j,k}  $. So $ \bar{g}_{p,j,k}(x) = 2^{j} \tilde{g}_{p,j,k} \chi_{j,k}(x)$. Now we put $ g_{p,j,k} =  ( p + 1 )^{(2m+3)^{2} }  2^{js} 2^{\frac{j}{2}}  \tilde{g}_{p,j,k}  $. It follows that $ g(\lambda)  $ can be written as in \eqref{seq_lin_fun}. Moreover we observe that we have the equality $ \bar{g}_{p,j,k}(x) =   g_{p,j,k}  ( p + 1 )^{- (2m+3)^{2} }  2^{-js} 2^{ \frac{j}{2}}   \chi_{j,k}(x)$. Consequently we get
\begin{align*}
\Vert g \vert  (  f^{s}_{r,q}(m)    )'   \Vert  = \Big\| \Big[ \sum_{p = 0}^{\infty}  \sum_{j= -1}^\infty \sum_{k \in \mathbb{Z}}     ( p + 1 )^{- q' (2m+3)^{2} }  2^{-jsq'} 2^{ \frac{jq'}{2}} | g_{p,j,k} |^{q'}  \chi_{j,k}(x)     \Big]^{\frac{1}{q'}}\Big|  L_{r'} (\mathbb{R}) \Big\| .
\end{align*}
This completes step 1 of the proof. 

\textit{Step 2.} Now we prove $   f^{-s}_{r',q'}(m) \subset (  f^{s}_{r,q}(m)    )'  $. For that purpose let $ \{ g_{p,j,k}  \} \subset \mathbb{C}    $ be a sequence with $ \{ g_{p,j,k}  \} \in  f^{-s}_{r',q'}(m)  $. We show that then \eqref{seq_lin_fun} is a bounded linear functional for all sequences $ \lambda = \{ c_{p,j,k}    \}  \in  f^{s}_{r,q}(m)   $. We observe
\begin{align*}
| g(\lambda) | &  \leq \int_{\mathbb{R}} \sum_{p = 0}^{\infty}  \sum_{j= -1}^\infty \sum_{k \in \mathbb{Z}} 2^{j}  | g_{p,j,k} |    | c_{p,j,k} |  \chi_{j,k}(x)   dx  \\
&  = \int_{\mathbb{R}} \sum_{p = 0}^{\infty}  \sum_{j= -1}^\infty \sum_{k \in \mathbb{Z}}  ( p + 1 )^{-  (2m+3)^{2} } 2^{-js} 2^{\frac{j}{2}}  | g_{p,j,k} |  ( p + 1 )^{ (2m+3)^{2} } 2^{js}  2^{\frac{j}{2}}   | c_{p,j,k} |  \chi_{j,k}(x)   dx  .
\end{align*}
Next we use the H\"older inequality with $ \frac{1}{q} + \frac{1}{q'} = 1   $. Then we find
\begin{align*}
| g(\lambda) |  &  \leq \int_{\mathbb{R}} \Big ( \sum_{p = 0}^{\infty}  \sum_{j= -1}^\infty \sum_{k \in \mathbb{Z}}  ( p + 1 )^{- q' (2m+3)^{2} } 2^{-jsq'} 2^{\frac{jq'}{2}}  | g_{p,j,k} |^{q'}  \chi_{j,k}(x)  \Big )^{\frac{1}{q'}}  \\
& \qquad \qquad \qquad \qquad \times \Big (   \sum_{p = 0}^{\infty}  \sum_{j= -1}^\infty \sum_{k \in \mathbb{Z}}    ( p + 1 )^{q (2m+3)^{2} } 2^{jsq}  2^{\frac{jq}{2}}   | c_{p,j,k} |^{q}  \chi_{j,k}(x)  \Big )^{\frac{1}{q}}   dx  .
\end{align*}
Now we apply the H\"older inequality again with $ \frac{1}{r} + \frac{1}{r'} = 1   $. We observe
\begin{align*}
| g(\lambda) |  &  \leq \Big \Vert \Big ( \sum_{p = 0}^{\infty}  \sum_{j= -1}^\infty \sum_{k \in \mathbb{Z}}  ( p + 1 )^{- q' (2m+3)^{2} } 2^{-jsq'} 2^{\frac{jq'}{2}}  | g_{p,j,k} |^{q'}  \chi_{j,k}(x)  \Big )^{\frac{1}{q'}} \Big \vert L_{r'} (\mathbb{R})  \Big \Vert \\
& \qquad \qquad \qquad \qquad \times \Big \Vert \Big (   \sum_{p = 0}^{\infty}  \sum_{j= -1}^\infty \sum_{k \in \mathbb{Z}}    ( p + 1 )^{q (2m+3)^{2} } 2^{jsq}  2^{\frac{jq}{2}}   | c_{p,j,k} |^{q}  \chi_{j,k}(x)  \Big )^{\frac{1}{q}} \Big \vert L_{r} (\mathbb{R})  \Big \Vert   .
\end{align*}
Recall that we have $ \{ g_{p,j,k}  \} \in  f^{-s}_{r',q'}(m)  $ and $ \{ c_{p,j,k}    \}  \in  f^{s}_{r,q}(m)   $. Therefore the right-hand side is finite. This completes the proof. 
\end{proof}

Now we are well prepared to prove the main result of this section. It reads as follows. 

\begin{prop}\label{mainresult_sneg}
Let $  s < 0  $, $  1 < r < \infty $, $ 1 < q < \infty  $ and $ m \in \mathbb{N}  $ with $ m \geq 2  $. Moreover we assume that the parameters fulfill one of the conditions $ (IV)  $ or $(V)$ that can be found in Theorem \ref{mainresult1}. Let $ f \in \mathcal{S}'(\mathbb{R})   $. Then we have   $ f \in F^{s}_{r,q}(\mathbb{R})  $ if and only if $ f $ can be represented as
\begin{equation}\label{rep_main_negs}
f = \sum_{p \geq 0} \sum_{j = -1}^{\infty} \sum_{k \in \mathbb{Z}} c_{p,j,k} \psi_{p,j,k}
\end{equation}
with convergence in $ \mathcal{S}'(\mathbb{R})  $, where we have
\begin{align*}
\Big\| \Big[ \sum_{p = 0}^{\infty}  \sum_{j= -1}^\infty \sum_{k \in \mathbb{Z}}  ( p + 1 )^{-(2m+3)^{2} q}  2^{jsq} 2^{\frac{jq}{2}}    |c_{p,j,k}|^{q}  |\chi_{j,k}(x)|^{q}   \Big]^{\frac{1}{q}}\Big|  L_{r} (\mathbb{R}) \Big\|   < \infty .
\end{align*}
Moreover the quasi-norms $ \Vert f \vert F^{s}_{r,q}(\mathbb{R})  \Vert    $ and
\begin{align*}
\inf_{\eqref{rep_main_negs}} \Big\| \Big[ \sum_{p = 0}^{\infty}  \sum_{j= -1}^\infty \sum_{k \in \mathbb{Z}}  ( p + 1 )^{-(2m+3)^{2} q}  2^{jsq} 2^{\frac{jq}{2}}    |c_{p,j,k}|^{q}  |\chi_{j,k}(x)|^{q}   \Big]^{\frac{1}{q}}\Big|  L_{r} (\mathbb{R}) \Big\|  
\end{align*}
are equivalent. Here the infimum is taken over all sequences $ \{ c_{p,j,k}  \}   $ such that \eqref{rep_main_negs} is fulfilled.
\end{prop}

\begin{proof}
\textit{Step 1. Preparations and parameters.} To prove this result we use the characterizations for positive smoothness in combination with some duality arguments. For that purpose we define $ 1 < r' < \infty $ and $ 1 < q' < \infty$ such that
\begin{equation}\label{par_dual_profneg}
\frac{1}{r} + \frac{1}{r'} = 1 \qquad \qquad \mbox{and} \qquad \qquad \frac{1}{q} + \frac{1}{q'} = 1 .
\end{equation}
Let $  s < 0  $. Then from Lemma \ref{Fsrq_dual} and Lemma \ref{fsrqm_dual} we know that we have 
\begin{align*}
(  F^{s}_{r,q}(\mathbb{R})    )' = F^{-s}_{r',q'}(\mathbb{R})  \qquad \qquad \mbox{and} \qquad \qquad (  f^{s}_{r,q}(m)    )' = f^{-s}_{r',q'}(m).
\end{align*}
Moreover from the Propositions \ref{res_Fspq_m-1}, \ref{res_Fspq_finalup} and \ref{Fspq_finallow} we know that there is a quarklet characterization similar to that given in Proposition \ref{mainresult_sneg} if we are in one of the following cases:
\begin{itemize}
\item[(A)] We have $0 < -s  < m - 1$ .

\item[(B)] We have $ -s \geq m - 1  $ with $  -s  < \min  ( m - 1 + \frac{1}{r'} , m - 1 + \frac{1}{q'}     )$. For the fine index $q'$ we assume $ \frac{1}{q'} < -s + 1 - mr'(-s + 1 - m)   $.
\end{itemize}
When we combine \eqref{par_dual_profneg} with (A) and (B) we obtain the conditions concerning the parameters that can be found in Proposition \ref{mainresult_sneg}.

\textit{Step 2. The upper estimate.} Now we prove that for $f$ given by \eqref{rep_main_negs} there exists a $ C > 0    $ independent of $f$ such that
\begin{align*}
& \|\, f \, |F^s_{r,q}(\mathbb{R})\| \leq C \Big\| \Big[ \sum_{p = 0}^{\infty}  \sum_{j= -1}^\infty \sum_{k \in \mathbb{Z}} ( p + 1 )^{-(2m+3)^{2} q}  2^{jsq} 2^{\frac{jq}{2}}    |c_{p,j,k}|^{q}  |\chi_{j,k}(x)|^{q}   \Big]^{\frac{1}{q}}\Big|  L_{r} (\mathbb{R}) \Big\| .
\end{align*}
For that purpose at first recall that we have $   F^{s}_{r,q}(\mathbb{R})     = ( F^{-s}_{r',q'}(\mathbb{R}) )'  $, see Lemma \ref{Fsrq_dual}. Therefore we can write $ \|\, f \, |F^s_{r,q}(\mathbb{R})\| = \sup \{ \vert f(g)   \vert : \|\, g \, |F^{-s}_{r',q'}(\mathbb{R})\|  \leq 1   \} $. That means we interpret $f$ as a linear functional. Notice that the parameters $ (-s,r',q')   $ are given in such a way that we can apply Proposition \ref{Fspq_finallow}. With other words we can find representations
\begin{equation}\label{rep_qua_prs<0}
g = \sum_{p \geq 0} \sum_{j = -1}^{\infty} \sum_{k \in \mathbb{Z}} g_{p,j,k} \psi_{p,j,k}
\end{equation}
and a constant $ C > 0  $ independent of $ g \in F^{-s}_{r',q'}(\mathbb{R})  $ such that
\begin{align*}
\inf_{\eqref{rep_qua_prs<0}} \Big\| \Big[ \sum_{p = 0}^{\infty}  \sum_{j= -1}^\infty \sum_{k \in \mathbb{Z}}  ( p + 1 )^{(2m+3)^{2} q'} 2^{-jsq'} 2^{\frac{jq'}{2}}    |g_{p,j,k}|^{q'}  |\chi_{j,k}(x)|^{q'}   \Big]^{\frac{1}{q'}}\Big|  L_{r'} (\mathbb{R}) \Big\|  \leq C  \|\, g \, |F^{-s}_{r',q'}(\mathbb{R})\| .
\end{align*}
Here the infimum is taken over all sequences $ \{ g_{p,j,k}  \}   $ such that \eqref{rep_qua_prs<0} is fulfilled. Consequently we also can write
\begin{align*}
& \|\, f \, |F^s_{r,q}(\mathbb{R})\| \\
& \leq \sup \Big \{ \vert f(g)   \vert : C_{1} \inf_{\eqref{rep_qua_prs<0}} \Big\| \Big[ \sum_{p,j,k}    ( p + 1 )^{(2m+3)^{2} q'} 2^{-jsq'} 2^{\frac{jq'}{2}}    |g_{p,j,k}|^{q'}  |\chi_{j,k}(x)|^{q'}   \Big]^{\frac{1}{q'}}\Big|  L_{r'} (\mathbb{R}) \Big\| \leq 1  \Big \} .
\end{align*}
Notice that we have sequence representations for both $f$ and $g$. Hence $  f(g)  $ also can be interpreted as a linear functional in the sense of \eqref{seq_lin_fun}. Next we use that we have $   f^{s}_{r,q}(m)    = ( f^{-s}_{r',q'}(m) )'   $, see Lemma \ref{fsrqm_dual}. So with Definition \ref{def_seqsp_qua} we get 
\begin{align*}
& \|\, f \, |F^s_{r,q}(\mathbb{R})\| \leq C_{2} \Big\| \Big[ \sum_{p = 0}^{\infty}  \sum_{j= -1}^\infty \sum_{k \in \mathbb{Z}} ( p + 1 )^{-(2m+3)^{2} q}  2^{jsq} 2^{\frac{jq}{2}}    |c_{p,j,k}|^{q}  |\chi_{j,k}(x)|^{q}   \Big]^{\frac{1}{q}}\Big|  L_{r} (\mathbb{R}) \Big\| .
\end{align*}
This completes step 2 of the proof.

\textit{Step 3. The lower estimate.} Let $  f \in F^s_{r,q}(\mathbb{R})   $. We will prove that there is a constant $ C > 0  $ independent of $f$ such that we have
\begin{align*}
\inf_{\eqref{rep_main_negs}} \Big\| \Big[ \sum_{p = 0}^{\infty}  \sum_{j= -1}^\infty \sum_{k \in \mathbb{Z}} ( p + 1 )^{-(2m+3)^{2} q}  2^{jsq} 2^{\frac{jq}{2}}    |c_{p,j,k}|^{q}  |\chi_{j,k}(x)|^{q}   \Big]^{\frac{1}{q}}\Big|  L_{r} (\mathbb{R}) \Big\| \leq C \|\, f \, |F^s_{r,q}(\mathbb{R})\| .
\end{align*}
Here the infimum is taken over all representations of the form \eqref{rep_main_negs}. To prove this we proceed in a similar way as it is described in the proofs of the Propositions \ref{thm_spline_trCW}, \ref{Fspq_prefinallow} and \ref{Fspq_finallow}. However because of the negative smoothness some modifications are necessary. Therefore below we will give some details.

\textit{Substep 3.1.} At first we prove that there is a representation 
\begin{equation}\label{eq-s<0-atoms1}
f = \sum_{j = -1}^{\infty} \sum_{k \in \mathbb{Z}} \lambda_{j,k} \psi_{0,j,k}
\end{equation}
and that we have the estimate 
\begin{equation}\label{eq-s<0-atoms2}
\inf_{\eqref{eq-s<0-atoms1}} \Big\| \Big[   \sum_{j= -1}^\infty \sum_{k \in \mathbb{Z}}   2^{jsq} 2^{\frac{jq}{2}}    |\lambda_{j,k}|^{q}  |\chi_{j,k}(x)|^{q}   \Big]^{\frac{1}{q}}\Big|  L_{r} (\mathbb{R}) \Big\| \leq C \|\, f \, |F^s_{r,q}(\mathbb{R})\| .
\end{equation}
Here the infimum is taken over all admissible representations \eqref{eq-s<0-atoms1}. To see this we want to interpret the functions $ \psi_{0,j,k}  $ as atoms. For a definition concerning atoms we refer to \cite{Tr10}, see Definition 1.5. We use the notation from there. It is not difficult to see that the family $ \psi_{0,j,k}  $ with $ j \geq -1  $ and $ k \in \mathbb{Z}  $ is consistent with that definition. So since we have $ \psi_{0,j,k} \in L_{\infty}(\mathbb{R})    $ we can put $ K = 0  $. Because of the moment conditions we can use $ L = m   $. Moreover since the functions $ \psi_{0,j,k}  $ are compactly supported there exists a $ d > 1   $ like it is described in Definition 1.5 in \cite{Tr10}. All in all we conclude that the functions $ \psi_{0,j,k}  $ can be seen as atoms. Therefore we can apply Theorem 1.7 from \cite{Tr10}. It tells us that we have a representation \eqref{eq-s<0-atoms1} and the estimate \eqref{eq-s<0-atoms2} in the case of $ -m < s < 0   $. But this condition is fulfilled. So substep 3.1 is complete.

\textit{Substep 3.2.} Now we prove the full result as we formulated it at the beginning of step 3. For that purpose we use the result from substep 3.1 in combination with the argument from the proof of Proposition \ref{Fspq_finallow}. Here only minor modifications are necessary. Therefore we omit the details. The proof is complete.    
\end{proof}

After all the preparations we performed so far, we are now in position to
prove the main result of this paper, namely Theorem \ref{mainresult1}.

\vspace{0,3 cm}

\textbf{Proof of Theorem \ref{mainresult1}.} This result is a combination of the Propositions \ref{res_Fspq_m-1}, \ref{res_Fspq_finalup}, \ref{Fspq_finallow} and \ref{mainresult_sneg}.

\begin{rem}\label{rem_nec1}
In the formulation of Theorem \ref{mainresult1} many conditions concerning the parameters show up. It is known that at least some of them are also necessary. So it is not difficult to see that the assertion from Theorem \ref{mainresult1} does not hold when we have $  s \geq m - 1 + \frac{1}{r}   $. This is a simple consequence of Lemma \ref{spline_inF}. For the case $ m = 1 $ and $ p = 0$, namely the case of the Haar system without polynomial enrichment, even more has been proved.  Here it is known that the Haar system is an unconditional basis in $   F^{s}_{r,q}(\mathbb{R})$ if and only if
\[
\left\{ \begin{array}{lll}
0 < r < \infty, 0 < q < \infty,  \qquad  \max (\frac{1}{r} , \frac{1}{q} , 1) - 1 < s < \min (\frac{1}{r} , \frac{1}{q} , 1) ;
\\
1 < r < \infty, 1 < q < \infty, \qquad s = 0 ;
\\  
1 < r < \infty, 1 < q < \infty, \qquad \max (\frac{1}{r} , \frac{1}{q} ) -1 < s < 0  .
\end{array}
\right.
\]
This result can be found in \cite{GSU}, see Theorem 1.1. A forerunner for the case $ 1 < r,q < \infty $ is given in \cite{SU17}, see Theorem 1.1. We also want to refer to Remark 3.21. in \cite{Tr20} where a short summary of the situation can be found. A corresponding result for the Besov spaces $   B^{s}_{r,q}(\mathbb{R})$ is formulated in \cite{Tr06}, see Theorem 1.58. In connection with that we also would like to refer to Theorem 1, Theorem 2 and Corollary 2 in \cite{Tr78Ser}.  
\end{rem}

\section{Quarklet Characterizations for Triebel-Lizorkin-Morrey spaces}\label{sec_mor}

In this section we want to prove quarklet characterizations for the more general Triebel-Lizorkin-Morrey spaces $ \mathcal{E}^{s}_{u,r,q}(\mathbb{R})   $. More precisely we will give a proof for Theorem \ref{mainresult2}. The spaces $ \mathcal{E}^{s}_{u,r,q}(\mathbb{R})    $ are function spaces that are built upon Morrey spaces. Because of this at first we want to recall the definition of the Morrey spaces $ \mathcal{M}^{u}_{r}(\mathbb{R})  $. 

\begin{defi}\label{def_mor}
Let $ 0 < r \leq u < \infty$. Then the  Morrey space $ \mathcal{M}^{u}_{r}(\mathbb{R})  $ is defined to be the set of all functions $ f \in L_{r}^{loc}(\mathbb{R}) $ such that 
\begin{align*}
\Vert f \vert \mathcal{M}^{u}_{r}(\mathbb{R}) \Vert := \sup_{- \infty \leq a < b \leq \infty} \vert b - a \vert^{\frac{1}{u}-\frac{1}{r}} \Big ( \int_{a}^{b} \vert f(x) \vert^{r} dx      \Big )^{\frac{1}{r}} < \infty.
\end{align*} 
\end{defi}
The Morrey spaces $ \mathcal{M}^{u}_{r}(\mathbb{R})  $ are quasi-Banach spaces and Banach spaces for $ r \geq 1$. They have many connections  to the Lebesgue spaces $ L_{r}(\mathbb{R})$. So for $ r \in (0,\infty) $ we have $ \mathcal{M}^{r}_{r}(\mathbb{R}) = L_{r}(\mathbb{R})$. Moreover for $ 0 < r_{2} \leq r_{1} \leq u < \infty $ we have $ L_{u}(\mathbb{R}) = \mathcal{M}^{u}_{u}(\mathbb{R}) \hookrightarrow   \mathcal{M}^{u}_{r_{1}}(\mathbb{R})  \hookrightarrow  \mathcal{M}^{u}_{r_{2}}(\mathbb{R}) $. Now let us define the one-dimensional Triebel-Lizorkin-Morrey spaces.
\begin{defi}\label{def_tlm}
Let $ 0 < r \leq u < \infty $, $ 0 < q \leq \infty $ and $ s \in \mathbb{R} $. $ (\lambda_{k})_{k\in \N_0 }$ is a smooth dyadic decomposition of the unity. Then the Triebel-Lizorkin-Morrey space $  \mathcal{E}^{s}_{u,r,q}(\mathbb{R}) $ is defined to be the set of all distributions $ f \in \mathcal{S}'(\mathbb{R})  $ such that
\begin{align*} 
\Vert f \vert \mathcal{E}^{s}_{u,r,q}(\mathbb{R})  \Vert := \Big  \Vert \Big ( \sum_{k = 0}^{\infty} 2^{ksq} \vert \mathcal{F}^{-1}[\lambda_{k} \mathcal{F}f](x) \vert ^{q} \Big   )^{\frac{1}{q}} \Big \vert \mathcal{M}^{u}_{r}(\mathbb{R})   \Big \Vert < \infty .
\end{align*}
In the case $ q = \infty $ the usual modifications are made.

\end{defi}
The Triebel-Lizorkin-Morrey spaces are generalizations of the original Triebel-Lizorkin spaces. So for $ u = r  $ we observe $ \mathcal{E}^{s}_{r,r,q}(\mathbb{R}) = F^{s}_{r,q}(\mathbb{R})  $. A first systematically collection of the properties of the spaces $ \mathcal{E}^{s}_{u,r,q}(\mathbb{R})  $ can be found in \cite{ysy}. Now let us prove Theorem \ref{mainresult2}. For that purpose in a first step we show the upper estimate. 

\begin{prop}\label{res_TLMs_R1}
Let $  0 < r \leq u < \infty $, $ 0 < q < \infty  $ and $ m \in \mathbb{N}  $ with $ m \geq 2  $. Let
\begin{equation}\label{s_in_Esurq_finallow}
\max \Big(0, \frac 1r -1 , \frac 1q -1  \Big) < s  < m - 1.
\end{equation}
Let $ f \in \mathcal{S}'(\mathbb{R})  $ such that there exists a representation 
\begin{equation}\label{eq-infTLM3}
f = \sum_{p \geq 0} \sum_{j \geq -1} \sum_{k \in \mathbb{Z}} c_{p,j,k} \psi_{p,j,k}
\end{equation}
with convergence in $ \mathcal{S}'(\mathbb{R})  $. Then there exists a $ C > 0  $ independent from $ f $ such that   
\begin{align*}
& \|\, f \, | \mathcal{E}^{s}_{u,r,q}(\mathbb{R})  \| \leq C \inf_{\eqref{eq-infTLM3}} \Big\| \Big[ \sum_{p = 0}^{\infty}  \sum_{j= -1}^\infty \sum_{k \in \mathbb{Z}} ( p + 1 )^{(2m+3)^{2} q}  2^{jsq} 2^{\frac{jq}{2}}    |c_{p,j,k}|^{q}  |\chi_{j,k}(x)|^{q}   \Big]^{\frac{1}{q}}\Big|  \mathcal{M}^{u}_{r}(\mathbb{R})  \Big\| .
\end{align*}
Here the infimum is taken over all sequences $ \{ c_{p,j,k}  \}_{p,j,k} \subset \mathbb{C} $ such that \eqref{eq-infTLM3} is fulfilled.
\end{prop}

\begin{proof}
This result can be proved in the same way as Proposition \ref{res_Fspq_m-1}. All arguments that are used there also hold for the more general Triebel-Lizorkin-Morrey spaces. We only have to make sure that all tools have counterparts for the spaces $ \mathcal{E}^{s}_{u,r,q}(\mathbb{R})   $. But this is the case. So the Triebel-Lizorkin-Morrey spaces can be described in terms of differences which means there is a counterpart for Lemma \ref{discrete}. For that we refer to \cite{HoTLM}, see Theorem 2 and Theorem 6. Here we require the condition $\sigma_{r,q} < s $. One might also consult \cite{ysy}, see Chapter 4.3.1. On the other hand we need the Hardy-Littlewood-Maximal inequality. Here we also have a counterpart for the Morrey case. For that we refer to \cite{TangXu}, see Lemma 2.5. All the other arguments used in the proof of Proposition \ref{res_Fspq_m-1} are independent from the question whether we work with the Lebesgue or the Morrey case. So we can proceed like there to obtain the desired result.  
\end{proof}

\begin{rem}
It seems that it is not possible to use complex interpolation like it is described in the proof of Proposition \ref{res_Fspq_finalup} to obtain a counterpart of Proposition \ref{res_Fspq_finalup} for the Triebel-Lizorkin-Morrey spaces. The reason for this is the more difficult behavior of complex interpolation in the context of the spaces $  \mathcal{E}^{s}_{u,r,q}(\mathbb{R})     $. More details concerning this topic can be found in \cite{zh_ho_si}, see in particular Proposition 1.6 in \cite{zh_ho_si}. 
\end{rem}

Now we prove the corresponding lower estimate.

\begin{prop}\label{Esurq_finallow}
Let $  0 < r \leq u < \infty $, $ 0 < q < \infty  $ and $ m \in \mathbb{N}  $ with $ m \geq 2  $. Moreover let \eqref{s_in_Esurq_finallow} be fulfilled.  Let $ f \in \mathcal{E}^{s}_{u,r,q}(\mathbb{R})   $.  Then there exists a sequence $ \{ c_{p,j,k}  \}_{p \in \mathbb{N}_{0}, k \in \mathbb{Z}, j \geq -1}   $ such that $ f $ can be represented as in \eqref{eq-infTLM3}. Moreover there is a constant $ C > 0  $ independent of $ f \in \mathcal{E}^{s}_{u,r,q}(\mathbb{R})   $ such that
\begin{align*}
\inf_{\eqref{eq-infTLM3}} \Big\| \Big[ \sum_{p = 0}^{\infty}  \sum_{j= -1}^\infty \sum_{k \in \mathbb{Z}}  ( p + 1 )^{(2m+3)^{2}q} 2^{jsq} 2^{\frac{jq}{2}}    |c_{p,j,k}|^{q}  |\chi_{j,k}(x)|^{q}   \Big]^{\frac{1}{q}}\Big| \mathcal{M}^{u}_{r}(\mathbb{R})  \Big\|  \leq C  \|\, f \, | \mathcal{E}^{s}_{u,r,q}(\mathbb{R})  \| .
\end{align*}
Here the infimum is taken over all sequences $ \{ c_{p,j,k}  \}   $ such that \eqref{eq-infTLM3} is fulfilled.
\end{prop}

\begin{proof}
This result can be proved in the same way as Proposition \ref{Fspq_finallow}. Therefore we will keep the proof rather short.

\textit{Step 1.} At first we prove a counterpart of Proposition \ref{thm_spline_trCW}. That means that for a function $ f \in \mathcal{E}^{s}_{u,r,q}(\mathbb{R})  $ we put $ \lambda_{j,k} = \langle   \tilde{\psi}_{j,k} , f  \rangle_{L_{2}(\mathbb{R})} $ and prove that there is a constant $ C_{1} > 0  $ independent of $ f \in \mathcal{E}^{s}_{u,r,q}(\mathbb{R})  $ such that
\begin{equation}\label{eq_localm_ros}
\Big \Vert \Big ( \sum_{j = -1}^{\infty} \sum_{k \in \mathbb{Z}} 2^{jsq} 2^{\frac{jq}{2}} |\lambda_{j,k}|^{q} |\chi_{j,k}(x)|^{q}    \Big )^{\frac{1}{q}} \Big \vert  \mathcal{M}^{u}_{r}(\mathbb{R})   \Big \Vert \leq C_{1}  \|\, f \, | \mathcal{E}^{s}_{u,r,q}(\mathbb{R})  \| .
\end{equation}
For that purpose we proceed like in the proof of Proposition \ref{thm_spline_trCW}. There we have seen that the functions $ 2^{\frac{j}{2}} \tilde{\psi}_{j,k} $ can be interpreted as kernels of local means. Fortunately the theory of local means also has been developed for the Triebel-Lizorkin-Morrey spaces. For that we refer to \cite{Ros}, see Section 3. So the functions $ 2^{\frac{j}{2}} \tilde{\psi}_{j,k}    $ are $[A,B,C]-$ kernels of local means in the sense of Definition 3.1 in \cite{Ros} when we put $ A = 0 $ and $ B = m     $. Since the functions $  \tilde{\psi}_{j,k}    $ are compactly supported we can choose $ C > 0 $ sufficiently large. Consequently we can use Theorem 3.10 in \cite{Ros}. So we obtain formula \eqref{eq_localm_ros} for $ \sigma_{r,q} < s < m $. This completes step 1.

\textit{Step 2.}
Now we prove a counterpart for Proposition \ref{Fspq_prefinallow} for all $s$ such that \eqref{s_in_Esurq_finallow} is fulfilled. With other words we prove that $ f \in \mathcal{E}^{s}_{u,r,q}(\mathbb{R})  $ can be represented as in \eqref{rep_ohnequa} such that \eqref{eq_localm_ros} is fulfilled. Here the constant $C_{1}$ is independent from $f$. To show this we can use the arguments that are described in the proof of Proposition \ref{Fspq_prefinallow} in combination with Proposition \ref{res_TLMs_R1}. All the required modifications in the proof are obvious. Therefore we omit the details. 

\textit{Step 3.} 
Now we prove the full result, namely Proposition \ref{Esurq_finallow} itself. For that purpose we follow the proof of Proposition \ref{Fspq_finallow} and use the results from Step 1 and Step 2. Notice that the arguments that are described in the proof of Proposition \ref{Fspq_finallow} are independent from the question whether we work with Triebel-Lizorkin or with Triebel-Lizorkin-Morrey spaces. Therefore we can proceed like there to obtain the desired result. The proof is complete.
\end{proof}

\textbf{Proof of Theorem \ref{mainresult2}.} This result is a combination of the Propositions \ref{res_TLMs_R1} and \ref{Esurq_finallow}.

\begin{rem}
There also exist characterizations in terms of biorthogonal wavelets for generalized Triebel-Lizorkin-type spaces that have been introduced in \cite{LSUYW}. Here we want to refer to Chapter 4.4 and especially to Theorem 4.12. Notice that the Triebel-Lizorkin-Morrey spaces $  \mathcal{E}^{s}_{u,r,q}(\mathbb{R}) $ fit into the theory described in \cite{LSUYW}, see Example 3.3 and Chapter 11.2 in \cite{LSUYW}.
\end{rem}

\vspace{0,3 cm}

\textbf{Funding.} This paper is a result of the DFG project "adaptive high-order quarklet frame methods for elliptic operator equations" with grant number $DA360/24-1$. The author Marc Hovemann is funded by this project.

\vspace{0,3 cm}

\textbf{Acknowledgment.} The authors would like to thank Thorsten Raasch, Dorian Vogel, Jonas Gadatsch and Winfried Sickel for several tips and hints.

\end{document}